\documentclass[UnnumRef]{book}

\usepackage{enghb2e}
\usepackage{amssymb}
\usepackage{amsmath}
\usepackage{graphics}
\usepackage{rotating}
\usepackage{color}

\def\S{Chapter }

\def\secnumwidth{23pt}
\makeatletter
\newif\ifchapdigdoub
\global\chapdigdoubfalse
\def\@makechapterhead#1{%
  \@chapdigit=\thechapter
  \@FMHead{\thechapter}{#1}%
  {\parindent\z@
    \def\author##1##2{\immediate\write\@tempfile{\string\author{##1}{##2}}}%
    \def\FolioBoldFont{}%
    \def\'{\noexpand\'}%
    \def\`{\noexpand\`}%
    \def\"{\noexpand\"}%
    \let\@b\bullet
    \def\bullet{\raisebox{2pt}{$\scriptscriptstyle\@b$}}%
    \let\SubsectionItalicFont\it
    \immediate\openout\@tempfile temp
    {\let\newline\par
      \let\chaptocbreak\par
      \def\hyperlink##1##2{##2}%
      \def\break{\string\break\ }%
      \global\setbox\@tempboxa\vbox{{%
        \leftskip23pt
\vspace*{-1.2pt}
\sloppy\baselineskip10pt\rightskip24pt plus 1fill\footnotesize%
        \CCTwoFont
        \makeatletter\@input{\thechapter.toc}\makeatother\par}}}%
    \immediate\closeout\@tempfile
    \@firstauthortrue
    \def\author##1##2{\if@firstauthor\@firstauthorfalse\else\vskip7pt\fi{\CAFont ##1}%
   \vskip3pt\def\@t{##2}\ifx\@t\@empty\else{\CAAFont ##2}\fi\vskip-2.5pt
}%
\parbox[b]{12pc}{\raggedright\baselineskip0pt\@input{temp}\ }
\box\@tempboxa}\vskip32pt}

\def\@mydottedtocline#1#2#3#4#5{
  \ifnum #1>\c@tocdepth
  \else
    \vskip 2.3pt plus.2\p@
    {\interlinepenalty\@M
      \leavevmode
      \def\@dotsep{1.2}%
      \@tempdima #3\relax
      \rightskip\z@
      \advance\hsize-\secnumwidth
\hskip-4pt
\vbox{\hsize231pt\ifnum#1=1\hskip-\secnumwidth\else\ifnum#1=3\else\fi\fi\CCOneFont \ifnum#1=2\hangindent\subsecnumwidth#4\else#4\fi\nobreak
        \if@pdf
        \else
          \leaders\hbox{{$\m@th\mkern\@dotsep mu.\mkern\@dotsep mu$}}\hfill\nobreak
          {\hbox to18\p@{\CCOneFont #5}}
\fi}\par}\fi}
\makeatother



\begin{document}

\title{Handbook of  Linear Algebra}


\def\ord#1{| #1 |} 
\def\comp #1{\overline{ #1 }} 

\def\R{\mathbb{R}}
\def\C{\mathbb{C}}
\def\Q{\mathbb{Q}}
\def\D{\mathbb{D}}
\def\Cn{\C^{n}}
\def\Qn{\Q^{n}}
\def\Qnn{\Q^{n\times n}}
\def\Cnn{\C^{n\times n}}
\def\Cmn{\complex^{m\times n}}
\def\Cmm{\complex^{m\times m}}
\def\Rnn{\R^{n\times n}}
\def\Fnn{F^{n\times n}}
\def\Fxnn{F[x]^{n\times n}}
\def\Fn{F^{n}}
\def\Rn{\R^{n}}
\def\JCF{Jordan canonical form}
\def\RCF{rational canonical form}
\def\RCFIF{invariant factors rational canonical form}
\def\RCFED{elementary divisors rational canonical form}
\def\tr{{\rm tr}}
\def\b{{\bf b}}
\def\x{{\bf x}}
\def\y{{\bf y}}
\def\z{{\bf z}}
\def\u{{\bf u}}
\def\w{{\bf w}}
\def\s{{\bf g}}
\def\0{{\bf 0}}
\def\df{\noindent{\bf Definitions:}}
\def\ap{\noindent{\bf Applications:}}
\def\fc{\noindent{\bf Facts:}}
\def\ex{\noindent{\bf Examples:}}
\def\L{{\lambda}}
\def\hL{{\hat \L}}
\def\Nnu{{N^{\nu}_{\lambda}}}
\def\J{{\EuScript J^\#}}


\def\mmtx#1{\renewcommand{\arraystretch}{1.0}\left[\begin{array}{cccccccccccc} #1
\end{array}\right]} 
\def\rank{{\rm rank}}
\def\sign{{\rm sign}}
\def\span{{\rm Span}}
\def\range{{\rm range}}
\def\ker{{\rm ker}}
\def\diag{{\rm diag}}
\def\LH{{\mathcal L}(H) }
\def\functionsR{{\mathcal R}}
\def\functionsC{{\mathcal C}}
\def\T{\mathbb{T}}
\def\N{\mathbb{N}}

\chapter[Spectral sets]{Spectral Sets}

\begin{chapterauthors}
\chapterauthor{Catalin Badea\vskip-1pt}{Universit\'e Lille 1}
\chapterauthor{Bernhard Beckermann\vskip-1pt}{Universit\'e Lille 1}
\end{chapterauthors}


\noindent 

\noindent Spectral sets and $K$-spectral sets,
introduced by John von Neumann in \cite{vN51}, offer a possibility to
estimate the norm of functions of matrices in terms of the sup-norm of the function.
Examples of such spectral sets include the numerical range or the pseudospectrum of a matrix,
discussed in Chapters 16 and 18. Estimating the norm of functions of matrices is an essential
task in numerous fields of pure and applied mathematics, such as (numerical)
linear algebra \cite{greenbaum,Hi08}, functional analysis \cite{Pa02}, and numerical analysis.
More specific examples include probability \cite{DelyonDelyon}, semi-groups and existence
results for operator-valued differential equations, the study of numerical
schemes for the time discretization of evolution equations \cite{Cr08}, or the convergence rate of GMRES (Section 41.7).

The notion of spectral sets involves many deep
connections between linear algebra, operator theory, approximation theory, and complex analysis.
One requires for example simple criteria such that the closed unit disk, the numerical range, or
another given set in the complex plane  is $K$-spectral for a given matrix $A$. In order to study
sharpness in matrix norm inequalities, one may look for extremal matrices or operators. It is also of
interest to consider joint spectral sets of several matrices. How are notions like functional calculus,
similarity, or dilation related to spectral sets? Is the intersection of spectral sets also spectral,
and what are optimal constants? How to apply the theory of spectral sets to the approximate computation
of matrix functions? The aim of this chapter is to give at least partial answers to these questions,
and to present a survey of the modern theory of spectral and $K$-spectral sets for operators and matrices.

\vspace*{1pc}
\section{Matrices and Operators}\label{s0}

Though many of the examples presented below are for matrices,
it is more natural to present the theory of spectral sets in terms of operators
acting on an abstract Hilbert space $H$. Those readers preferring matrices can always think of the
Hilbert space $\mathbb C^n$ equipped with the usual scalar product.

For the convenience of the reader, some few basic properties of operator theory are collected below.
Standard references are the books \cite{RiNa,FoiasSzNagy,NikolskiEasy}. All definitions reduce to known ones in the case of matrices.

\smallskip

\df

\smallskip

\noindent

For a complex Hilbert space $H$ with scalar product $\langle \cdot,\cdot\rangle$ we denote by $\LH$ the normed space of all bounded linear operators on $H$. The {\bf operator norm} of $A\in \LH$ is defined by $\|A\| = \sup \{ \|Ax\| : \|x\| = 1 \}$. A {\bf contraction} is a linear operator whose operator norm is not greater than one, while a {\bf strict contraction} is an operator $A$ such that $\|A\| < 1$.
An operator $A \in \LH$ is {\bf normal} if $A$ commutes with
its Hilbert space adjoint  $A^\ast$. The operator $A$ is {\bf unitary} if $AA^\ast =
A^\ast A = I$, where $I$ represents the identity operator. The operator $A$ is {\bf Hermitian},
or {\bf self-adjoint} if $A^\ast = A$. We say that the operator $A \in \LH$ is {\bf positive semidefinite}, and we write $A \succeq 0$, if $\langle Ax,x\rangle \ge 0$ for each $x\in H$.

The {\bf spectrum} $\sigma(A)$ of $A$ is the set
of all complex numbers $z$ such that $A - zI$ is not an invertible operator. The {\bf spectral radius} $\rho(A)$ of the
operator $A$ is given by $\rho(A) = \sup\{ |\lambda| : \lambda \in \sigma(A)\} $. The {\bf numerical range}, or field of values,
of an operator $A$ is defined
by $W(A) :=  \{\langle Ax,x\rangle : \|x\| = 1\}$. The {\bf numerical radius} of $A \in \LH$ is defined
by $w(A) :=  \sup\{ |\lambda| : \lambda \in W(A)\} $. An operator $A$ is a {\bf numerical radius contraction} if $w(A) \le 1$.

For $s > 0$, we say that the operator $A \in \LH$ belongs to the class $C_{s}$ of Sz.-Nagy and Foias if the inequality
$ \frac{2-s}{s} \|\zeta Ax\|^2 + 2\frac{s - 1}{s} \mathrm{re}( \langle \zeta Ax, x\rangle) \le \|x\|^2$
holds true for every $x\in H$ and every complex number $\zeta$ with $|\zeta| < 1$.
The {\bf operator radius} $w_{s}$ associated with the class
$C_{s}$ may be defined by
$w_{s}(A) = \inf \left\{ r : r> 0, \frac{1}{r}A \in C_{s}\right\} .$

\smallskip

The following notation will be used throughout this chapter: for a set $M \subset \mathbb C$ we denote by $\partial M$ its boundary and by $\overline{M}$ its closure.
$\mathbb D$ is the open unit disk, $\T = \partial \D$ is the unit circle and
$X$ always denotes a closed subset of $\mathbb C$.  

\medskip

\df

\smallskip

\noindent

$\functionsR (X)$ and $\functionsC (X)$ are the sets of complex-valued {\bf bounded rational functions} on $X$, and complex-valued {\bf bounded continuous functions} on $X$, respectively, equipped with the {\bf supremum norm} $\|f\|_X = \sup \{ |f(x)| : x \in X\}$.


\bigskip

\fc

\noindent
  The following facts can be found in \cite{FoiasNagySzeged, FoiasSzNagy}.
\begin{enumerate}
\item  We have $ \|A\|/2 \le w(A) \le \|A\|$, and $\|A\|/s \le w_{s}(A)$. See also Chapter 18.

\item $w_{s}$ is a norm on $\LH$ if and only if $0 < s \le 2$.

\item $\lim_{s\to\infty} w_{s}(A) = \rho(A)$ for every $A\in \LH$.

\item The class $C_1$ is the class of all contractions and $w_1(A) = \|A\|$.

\item The class $C_2$ is the class of all numerical radius contractions and $w_{2}(A) = w(A)$.
\end{enumerate}

\medskip

\ex\vspace{-5pt}

\begin{enumerate}
\item The
    space $\ell^2$ of complex square summable sequences $x=(x_0,x_1,x_2,\ldots)$ with scalar product $\langle x,y\rangle = \sum_{j=0}^\infty x_j\overline{y_j}$ and norm $\| x \| = \sqrt{\langle x,x\rangle}$ is an example of an infinite dimensional Hilbert space. Here the action of a linear operator can be described by the matrix product $Ax$ with an infinite matrix $A$. The shift operator $S$ acting on $\ell^2$ via $S(z_0, z_1, \ldots) = (0, z_0, z_1, \ldots)$ has the norm given by $\| S \| = 1$, the spectrum $\sigma(S)=\overline{\mathbb D}$, and the numerical range $W(S)=\mathbb D$.

\item We will sometimes also make use of Banach spaces
    like the space $\ell^p$ of complex sequences $x=(x_0,x_1,x_2,\ldots)$ with norm $\| x \| = \bigl( \sum_{j=0}^\infty |x_j|^p \bigr)^{1/p}$ where again the action of bounded linear operators is described by infinite matrices and the
    matrix-vector product. If one considers the subspaces of sequences where only the first $n$ entries are non-zero, then there is a canonical isomorphism with $\mathbb C^n$, equipped with the vector H\"older $p$-norm. Similarly, if the infinite matrix $A$ contains  non-zero entries only in the first $n$ rows and columns, then the operator norm $\| A \|_{\mathcal L(\ell^p)}$ coincides with the matrix H\"older $p$-norm of the corresponding principal submatrix of order $n$.

\item
    Another example of an infinite dimensional Banach space is the set $L^p$ of complex-valued functions defined on some set $X$, with $\| f \| := \bigl( \int_X |f(z)|^p \, d\mu(z) \bigr)^{1/p}$ for a suitable measure $\mu$ on $X$.

\item Since $f\in \functionsR (X)$ is analytic on a neighborhood of $X$, we get from the maximum principle for analytic functions that
    $$\|f\|_X = \sup \{ |f(x)| : x \in X\}  = \sup \{ |f(x)| : x \in
    \partial X\}.$$
    In particular, $\functionsR (X)$ can be seen as
    a subalgebra of $\functionsC (\partial X)$.
\end{enumerate}

\vspace*{1pc}
\section{Basic Properties of Spectral Sets}\label{s1}

Let $H$ be a complex Hilbert space and suppose that $A$ is a bounded
linear operator on $H$. Suppose now
that $\sigma(A)$ is included in the closed set $X$ and that $f=p/q \in \functionsR (X)$. As the poles of the rational
 function $f$ are outside of $X$, the operator $f(A)$ is
 naturally defined as $f(A) = p(A)q(A)^{-1}$
 (see Chapter 11)
 or, equivalently,
 by the Riesz holomorphic functional
 calculus \cite{RiNa}.

We can now define the central objects of study of this chapter.

\bigskip

\df

\smallskip

\noindent Let $H$ be a complex Hilbert space and suppose that $A \in \LH$ is a bounded linear operator on $H$.
Let $X$ be a closed set in the complex plane.

For a fixed constant $K > 0$, the set $X$ is said to be a
 $K$-{\bf spectral} set for $A$ if the spectrum $\sigma(A)$ of $A$
 is included in $X$ and the inequality $\|f(A)\| \le K\|f\|_X$ holds
 for every $f \in \functionsR (X)$.

The set $X$ is a {\bf spectral} set
 for $A$ if it is a $K$-spectral set with $K=1$.


\bigskip

\fc

\noindent Facts requiring proof for which no specific reference is given can be found in \cite[Chapter XI]{RiNa}
or \cite{vN51}.

\smallskip

 \begin{enumerate}

\item\label{c107s1ex1} Von Neumann inequalities for closed disks of the Riemann sphere:
  \begin{enumerate}
\item \label{c107s1ex1a}  a
closed disk $\{z \in \C : |z-\alpha| \le r\}$ is a spectral set
for $A\in \LH$ if and only if $\|A -\alpha I \| \le r$.

\item the closed set $\{z \in \C : |z-\alpha| \ge r\}$ is spectral for
$A \in \LH$ if and only if $\|(A -\alpha I)^{-1} \| \le
r^{-1}$.

\item \label{c107s1ex1c} the closed right half-plane $\C_0^+ = \{ \mathrm{re}(z)\geq 0 \}$ is a spectral set
for $A$ if and only if $\mathrm{re}(\langle A v,v\rangle)\geq0$ for all $v\in H$.
More generally, any closed halfplane is spectral for $A$ if and only if it contains the numerical range $W(A)$.
\end{enumerate}


\item \label{c107s1fa1}
Suppose that $X$ is a spectral set for $A$ and let
  $f \in \functionsR (\sigma(A))$.
Then $f(X)$ is spectral for $f(A)$.
More generally, if $(f_n) \in \functionsR(X)$, $\lim\|f_n - f\|_{X} = 0$, and $\lim\|f_n(A) - B\| = 0$,
then $f(X)$ is a spectral set for $B$.

\item \label{c107s1fac3} Any closed superset of a spectral set is again spectral. Any closed superset of a $K$-spectral set is again $K$-spectral.

\item
  The spectrum $\sigma(A)$ of $A\in \LH$ is the intersection of all spectral sets for $A$.

\item \cite{WilliamsSzeged} Any spectral set contains a \emph{minimal spectral set}, i.e., a spectral set having no proper closed subset which is spectral.

\item \cite{vN51} If the operator $A \in \LH$ is normal, then
the spectrum $\sigma(A)$
is a (minimal) spectral set for $A$.

\item \cite{WilliamsSzeged} Let $A$ be a Hilbert space operator, $X$ a set
containing $\sigma(A)$, and let $z_0$ be an interior point of $X$. If $\|f(A)\| \le \|f\|_X$ for
each rational function $f\in \functionsR(X)$ which vanishes at $z_0$, then $X$ is a spectral set for $A$.

\item \cite[page 18]{Pa02} Let $A$ be a Hilbert space operator and let $X$ be a closed set in the complex plane.
Let $\mathcal{S} :=  \functionsR (X) + \overline{\functionsR (X)}$ regarded as a subset of $\functionsC (\partial X)$.
If $X$ is spectral for $A$, then the
functional calculus homomorphism from $\functionsR (X)$ to $\mathcal{L}(H)$ defined by $f \mapsto f(A)$ extends to a well defined, \emph{positive map} (i.e., it sends positive functions to
positive semidefinite operators) $\Phi$ on $\mathcal{S}$ which sends $f+\overline{g}$ to $f(A)+g(A)^{\ast}$.
 Conversely, if $\Phi$ is well defined on $\mathcal{S}$, sends $f+\overline{g}$ to $f(A)+g(A)^{\ast}$ and is a positive map, then $X$ is spectral for $A$.

 \item Spectral sets and matrices with structure :
\begin{enumerate}
\item \label{c107s1fac9a}  A complex matrix $A\in \mathbb C^{n\times n}$ is normal if and only if
its set of eigenvalues $\sigma(A)$ is a spectral set for $A$.

\item \label{c107s1fac9b}  The unit circle $\T$ is a spectral set for $A\in \LH$ if and only
if $A$ is unitary.

\item \label{c107s1fac9c}  The real axis $\R$ is a spectral set for $A\in \LH$ if and only if $A$ is     self-adjoint 
    (= Hermitian).

\item \cite{WilliamsSzeged} Let $A$ be a completely non normal matrix, that is, the triangular
factor in a Schur decomposition of $A$ is not block diagonal. If $\| A \|=1$ then the
closed unit disk is a minimal spectral set.
\end{enumerate}

\item \label{c107s2faCrho} Let
   $s > 0$.
   If $A$ belongs to the class $C_{s}$, then the closed unit disk is $K$-spectral for $A$. \cite{FoiasNagySzeged}
We can take $K = 2s - 1$ if $s \ge 1$. \cite{OkuboAndo} $K=\max (1,s)$
is the best possible constant.

\item \cite{NevanlinnaJFA}  Lemniscates as $K$-spectral sets:
let $p$ be a polynomial with distinct roots and let $A\in \LH$. Let $R\ge 0$ satisfy $\|p(A)\|
\le R$ and be such that the lemniscate $\{z\in \C : |p(z)| = R\}$ contains no critical points of $p$. Then
$\{z\in \C : |p(z)| \le R\}$ is a $K$-spectral set for $A$.

\item \cite{KatsMats} The closed disk $3\overline{\D}$ of radius $3$ is spectral for
every Banach space contraction: if $A$ is an operator acting on a complex Banach
space $E$ such that $\|A\|_{\mathcal L(E)} \le 1$, and $p$ is a polynomial,
then $\|p(A)\|_{\mathcal L(E)} \le \sup\{|p(z)| : |z| \le 3\} $. The constant $3$ is the best possible one.
 \end{enumerate}

\bigskip

\ex\vspace{-5pt}


\noindent
\begin{enumerate}
\item \label{c107s1exa1}
   Suppose that $-1$ is not in the spectrum of the square matrix $A$ of order $n$. If the right half-plane $\C_0^+ = \{ \mathrm{re}(z) \geq 0 \}$ is a spectral set for $A$ then $\|f(A) \| \le \| f \|_{\C_0^+}=1$ for $f(z)=\frac{z-1}{z+1}$ by definition of a spectral set. In order to see that also the converse is true, suppose that $\| f(A) \| \leq 1$. We set $\mathbf{u}=(A+I)\mathbf{v}$, and observe that
   $$
     0 \leq \| \mathbf{u} \|^2 - \| f(A) \mathbf{u} \|^2 = \| (A+I)\mathbf{v} \|^2 - \| (A-I)\mathbf{v} \|^2 =
     4 \, \mathrm{re}(\langle A\mathbf{v},\mathbf{v}\rangle) .
   $$
   In particular, $\C_0^+$ contains the spectrum of $A$, and $\C_0^+$ is spectral for $A$ by Fact \ref{s1}.\ref{c107s1ex1c}.

\item
   Let $A \in \mathcal L(H)$, and
   $\alpha,\beta,\gamma,\delta\in \mathbb C$, $\alpha\delta-\beta\gamma\neq 0$, $-\delta/\gamma \not\in \sigma(A)$. Then both
   $f(z)=\frac{\alpha z + \beta}{\gamma z+\delta}$ and its inverse are rational functions. Applying twice Fact \ref{s1}.\ref{c107s1fa1}, we see that $X$ is spectral for $A$ if and only if $f(X)$ is spectral for $f(A)$.
   Thus for a proof of Fact \ref{s1}.\ref{c107s1ex1} one only needs to show that $\| A \|\leq 1$ implies that $\overline{\D}$ is spectral for $A$, see for instance \cite[Section 153]{RiNa}.

\item \label{c107s1exa3} Let the matrix $A$ be block diagonal with $A=\diag(A_1,A_2)$. Then we have $f(A)=\diag(f(A_1),f(A_2))$ and hence $X$ is $K$-spectral for $A$ if and only if it is $K$-spectral for $A_1$ and $A_2$.

\item \label{c107s1exa4} If the matrix $A$ can be factorized as $A=C B C^{-1}$, then $f(A)=C f(B) C^{-1}$.
    In particular, if $X$ is $K$-spectral for $B$ then it is also $K'$-spectral for $A$ with $K'\leq K \, \| C \| \, \| C^{-1} \|$.

\item \label{c107s1exa5} If the matrix $A$ is not diagonalisable, then no finite set $X$ can be $K$-spectral for $A$. For a proof, according to Example~\ref{s1}.\ref{c107s1exa3} and Example~\ref{s1}.\ref{c107s1exa4} we may suppose that $A$ is a Jordan block of order $\geq 2$. Suppose that $X=\{ x_1,...,x_k \}$ is finite.
    Direct computations show that $f(A)\neq 0$ for $f(x)=(x-x_1) \cdots (x-x_k)$, and hence $X$ cannot be $K$-spectral for $A$.

    It follows that the spectrum $\sigma(A)$ of a matrix $A$ is $K$-spectral if and only if $A$ is diagonalisable, and in this case we can take as $K$ the condition number of the matrix of eigenvectors. See also Example~\ref{sDilat}.\ref{diagonalizable}.

\item From Example~\ref{s1}.\ref{c107s1exa4} we see that a normal, Hermitian, or unitary matrix $A$ has $\sigma(A),\mathbb R$, or $\mathbb T$, respectively, as a spectral set. Let us show the converse result claimed in Facts~\ref{s1}.\ref{c107s1fac9a}--~\ref{s1}.\ref{c107s1fac9c}.

    If $\mathbb R$ is spectral for $A$, then according to Facts~\ref{s1}.\ref{c107s1fa1} and \ref{s1}.\ref{c107s1fac3}  the right half plane $\mathbb C_0^+$ is spectral for $iA$ and for $-iA$, which together with Fact~\ref{s1}.\ref{c107s1ex1c} implies that $\langle A \mathbf{u}, \mathbf{u} \rangle \in \mathbb R$ for all $\mathbf{u}$. Hence $A$ is Hermitian.

    If $\mathbb T$ is spectral for $A$, then $f(\mathbb T)=\mathbb R$ is spectral for $f(A)$, $f(z)=\frac{1}{i}\frac{z-1}{z+1}$, according to Fact~\ref{s1}.\ref{c107s1fa1}. Hence $f(A)$ is Hermitian, implying that $A$ is unitary.

    Finally, suppose that $\sigma(A)=\{\lambda_1,...,\lambda_n\}$ is spectral for $A$. Then, denoting by $\ell_1,...,\ell_n$ the Lagrange interpolation polynomials for $\sigma(A)$, that is, $\ell_j$ a polynomial of degree $\leq n-1$ with $\ell_j(\lambda_k)=0$ for $j\neq k$ and $=1$ for $j=k$, we find that
    $$
         z = \sum_{j=1}^n \lambda_j \ell_j(z) \quad \mbox{and hence} \quad
         A = \sum_{j=1}^n \lambda_j \ell_j(A).
    $$
    By Fact~\ref{s1}.\ref{c107s1fa1}, $\ell_j(\sigma(A))=\{ 0,1 \}\subset \mathbb R$ is spectral for $\ell_j(A)$, implying that $\ell_j(A)$ is Hermitian. From Example~\ref{s1}.\ref{c107s1exa5} we know that $A$ is diagonalisable, implying that the $\ell_j(A)$ commute, and thus $A$ is normal.

\item The closure $X$ of the $\epsilon$-pseudospectrum of $A\in \LH$ is a $K$-spectral set for
$K=\frac{\mbox{\small length}(\partial X)}{2\pi\epsilon}$ for any $\epsilon>0$, see Fact 16.3.5 in
Chapter 16.

 \item \label{c107s5fa1}

     \cite{BergerStampfli, OkuboAndo, Cr04a}
    The disk $\{ z\in \mathbb C: |z|\leq w(A) \}$ is $2$-spectral for $A$.
    $K=2$ is the best possible constant. This follows from Fact \ref{s1}.\ref{c107s2faCrho}.

\end{enumerate}

\vspace*{1pc}
\section{Around the von Neumann Inequality}\label{s1bis}
According to Fact~\ref{s1}.\ref{c107s1ex1a} stated for $\alpha = 0$, if $A$ is a Hilbert space contraction
and $f$ is a rational function with poles off the closed unit disk, then $\|f(A)\| \le \|f\|_{\overline{\D}}$. Since polynomials
are dense in the \emph{disk algebra} (that is, the Banach algebra of all
complex-valued functions which are analytic on $\D$ and continuous up to $\partial \D$),
the previous inequality implies that for every contraction on Hilbert space the rational
functional calculus extends to a functional calculus on the disk algebra.
The inequality $\|f(A)\| \le \|f\|_{\overline{\D}}$
is known in Operator Theory as the {\em von Neumann inequality}. The aim of this section is to present several variations and generalizations.

\smallskip

\fc

\smallskip

\begin{enumerate}

  \item \cite{foias57} The von Neumann inequality characterizes Hilbert spaces:
   if $E$ is a complex Banach space such that $\|p(A)\|_{\mathcal L(E)} \le \|p\|_{\overline{\D}}$ holds
 for every polynomial $p$ and every $A \in \mathcal{L}(E)$ with $\|A\|_{\mathcal L(E)} = 1$,
then $E$ is isometrically isomorphic to
a Hilbert space (i.e. the norm of $E$ comes from an inner product).

\item \cite{rovnyak82} Another converse of von Neumann inequality: suppose that $f_0$ is holomorphic on an open subset $G$ of $\D$ and $\|f_0(A)\|\le 1$ for every
contraction $A$ on a Hilbert space, with spectrum contained in $G$ (the operator $f_0(A)$ being defined by the Riesz
holomorphic functional calculus). Then $f_0$ is the restriction to $G$ of a holomorphic function $f$
defined and bounded by $1$ on $\D$.

\item \cite{KatsMats} The von Neumann inequality for arbitrary matrix norms:
let $p$ be a polynomial of degree $d$ with complex coefficients and let $A$ be a complex matrix
in $\mathbb C^{n\times n}$ with $\| A \|\leq 1$ for some subordinate matrix norm $\|\cdot\|$.
Then $\|p(A)\| \le (\pi n + 1)\, \|\, p\, \|_{\overline{\D}}$.

\item  The shift as an extremal operator: define the shift $S$ on a suitable space of
sequences by $S(z_0, z_1, \ldots) = (0, z_0, z_1, \ldots)$.

\begin{enumerate}
\item  \cite{HBohrineq, DixonvN} An operator-theoretical interpretation of H. Bohr's inequality:
let $r > 0$, let $E$ be a Banach space and let $A\in \mathcal{L}(E)$ be such
that $\|A\|_{\mathcal L(E)} \le r$.  Then $\|p(A)\|_{\mathcal L(E)} \le \|p(rS)\|_{\mathcal L(\ell^1)}$ for
every polynomial $p$.
We have $\|p(rS)\|_{\mathcal L(\ell^1)} \le \|p\|_{\overline{\D}}$ if and only if $r\le 1/3$.

\item \cite{Peller1978,CoiffRochWeiss78, NikolskiEasy} Matsaev inequality for some classes of contractions:
let $p$ be a real number between $1$ and $\infty$.
Let $A : L^p \mapsto L^p$ be an isometry on a $L^p$ space.
Then $A$ verifies the Matsaev inequality $\|f(A)\|_{\mathcal L(L^p)} \le \|f(S)\|_{\mathcal L(\ell^p)}$
for every polynomial $f$. This inequality also holds for
contractions on $L^p$ which preserve positive functions, or
for \emph{disjoint contractions} ($A(f)A(g) = 0$ whenever $fg = 0$), or for operators
such that
$\|A\|_{\mathcal L(L^1)} \le 1$ and $\|A\|_{\mathcal L(L^{\infty})} \le 1$.

\item \cite{Drury} Counterexample to a conjecture of Matsaev:
let $f(z) = 1 + 2z - \frac{22}{5}z^2$. There is a $2\times 2$ real matrix $A$
with a real (or complex) H\"older $4$-norm bounded above by $1$
but with the H\"older $4$-norm of $f(A)$ exceeding $\| f(S) \|_{\mathcal L(\ell^4)}$.
\end{enumerate}

\item Constrained von Neumann inequalities
\begin{enumerate}
\item \cite{PtakYoung}
Let $f$ and $g$ be two polynomials. Suppose that $A$ is a Hilbert space contraction with spectrum included in $\D$
such that $g(A) = 0$. Then $\|f(A)\| \le \|f(S^{\ast}\mid \ker \, g(S^{\ast}))\|$,
where $S^{\ast}$ is the backward shift,
$$ S^{\ast}(z_0, z_1, z_2, \ldots) = (z_1, z_2, \ldots),$$
which is
the adjoint of $S$ acting on $\ell^2$.

\item \cite{HaagerupdelaHarpe} Let $A$ be a Hilbert space nilpotent contraction with $A^n = 0$, $n \ge 2$. Then
$$w(A) \le
w(S^{\ast}\mid \ker \, (S^{n\ast}))) = \cos \frac{\pi}{n+1} .$$

\item \cite{BadeaCassier} Let $n\ge 2$. Let $A \in \LH$ be a contraction such that $A^n = 0$. Then
for each $s > 0$ and each polynomial $f$, we have $w_{s}(f(A)) \le w_{s}(f(S_n^{\ast}))$. Here
$S_n^{\ast}$ denotes the nilpotent Jordan block
$$S_n^{\ast} = \left[\begin{array}{r@{\quad}r@{\quad}r@{\quad}r@{\quad}r}
0 & 1 & 0 & \cdots & 0 \cr
			0 & 0 & 1 & \ddots & \vdots \cr
			\vdots & & \ddots & \ddots &  0\\
			\vdots &  &  & \ddots & 1 \\
			0 & \cdots & \cdots& \cdots & 0\end{array}\right] $$
which is unitarily equivalent to $S^{\ast}\mid \C^n = S^{\ast}\mid \ker \,S^{n\ast}$. In particular, for any $m$ we have
$$w(A^m) \leq \cos \frac{\pi}{k(m,n)+2}, \quad k(m,n) :=
[\frac{n-1}{m}].$$
Also, if the degree of the polynomial $f$ is at most $n-1$, then
$$
w_{s}(f(A))\leq \left(\frac{2}{s }-1\right)\left\| f\right\| _{\overline{\mathbb D}
}^{s}\left[ \inf_{\theta \in \R}\sup \{\left| f(\zeta )\right| :\zeta \in
\C,\zeta ^{2n-1}=e^{i\theta }\}\right] ^{1-s }$$
if $s \in (0,1]$, and
$$
w_{s }(f(A))\leq \left\| f\right\| _{\overline{\mathbb D}}^{2-s }\left[
\inf_{\theta \in \R}\sup \{\left| f(\zeta )\right| :\zeta \in \C,\zeta
^{2n-1}=e^{i\theta }\}\right] ^{s -1}$$
if $s \in [1,2]$.

\item \cite{BadeaCassier} Let $A\in \LH$ be a nilpotent Hilbert space operator satisfying
$$
  I - A^*A \succeq 0 , \quad I - 2A^*A + A^{*2}A^2 \succeq 0,
$$
and $A^n = 0$, $n
\geq 2$.
Then
$$w_{s}(f(A)) \leq w_{s}(f(B_{n}^*))$$
for all $s > 0$ and all polynomials $f$. Here $B_{n}^*$ is
given by the matrix
$$B_{n}^* = \left[\begin{array}{r@{\quad}r@{\quad}r@{\quad}r@{\quad}r}0 & \sqrt{\frac{1}{2}} & 0 & \cdots & 0 \\
			0 & 0 & \sqrt{\frac{2}{3}} & \ddots & \vdots \\
			\vdots & & \ddots & \ddots &  0\\
			\vdots &  &  & \ddots & \sqrt{\frac{n}{n+1}} \\
			0 & \cdots & \cdots& \cdots & 0\end{array}\right]
$$
which is unitarily equivalent to a compression of the Bergman shift.

\item \cite{BadeaCassier} Suppose $A \in \LH$ satisfies $\|A\| \leq 1$, $A^3 = 0$ and
$I - 2A^*A + A^{*2}A^2 \succeq  0$.
Then
$$w(A) \leq \sqrt{\frac{7}{24}} \quad \mbox{ and } \quad
w(A^2) \leq \sqrt{\frac{1}{12}}$$
and these constants are the best possible ones.

\item \cite{BadeaCassier} The link between the extremal operator in the constrained von Neumann
inequalities and the Taylor coefficients of rational functions positive
on the unit circle $\T$:
let $F=P/Q$ be a rational function
with no principal part and which is positive on the unit circle. Then the $k$th Taylor
coefficient $c_{k}$ satisfies
$
|c_{k}| \leq c_{0}w(R^{k}),
$
where $R=S^{\ast }\mid \ker (Q(S^{\ast }))$. In particular, if $P(e^{it})=\sum_{j=-n+1}^{n-1}c_{j}e^{ijt}$
is a positive trigonometric polynomial, $n\geq 2$, then
$
\left| c_{k}\right| \leq c_{0}\cos (\frac{\pi }{[\frac{n-1}{k}]+2})$
for $1\leq k\leq n-1
$
and, for every distinct numbers $k$ and $l$ among $\{0,\dots,n-1\}$, there
exists $\gamma \in \R$ such that
$
\left| c_{k}\right| +\left| c_{l}\right| \leq c_{0}
w(S_{n}^{k}+e^{i\gamma }S_{n}^{l}).
$
We also have
$$
\left| c_{k}\right| +\left| c_{l}\right| \leq c_{0}\left(1+\cos \frac{\pi }{[
\frac{n-1}{k+l}]+2}\right)^{1/2}
\left(1+\cos \frac{\pi }{[\frac{n-1}{\left| k-l\right|
}]+2}\right)^{1/2}.
$$
\end{enumerate}
%
%
%
\end{enumerate}

\ex
\smallskip

\begin{enumerate}

\item \cite[page 436]{RiNa}
The von Neumann inequality fails for the H\"older $1$-norm in $\mathbb C^2$ (or in $\ell^1$): consider
$$A= \left[\begin{array}{r@{\quad}r}
   0 & 1 \\
   1 & 0 \end{array}\right], \quad f(z) = \frac{z+a}{1+\overline{a}z}, \quad a=\frac{i}{2},
   \quad f(A) = \left[\begin{array}{r@{\quad}r}
   4i/5 & 3/5 \\
   3/5 & 4i/5 \end{array}\right],
   $$
then $\| A \|_1=1$, $\| f \|_{\overline{\mathbb D}}=1$, but $\| f(A) \|_1 =\frac{7}{5}>1$.

\item
\cite[page 23]{Pa02} The Schwarz-Pick lemma as a consequence of von Neumann inequality:
Let $p$ be a polynomial such that $\|p\|_{\overline{\D}} < 1$.
Let $a,b$ and $c$ be complex numbers such that $|a| < 1, |c| < 1$, $a\neq c$, and $|b|^2 = (1-|a|^2)(1-|c|^2)$.
Consider the $2\times 2$ matrix $$A= \left[\begin{array}{c@{\quad}c}
a & ~b \\
0 & ~c \end{array}\right]
$$ acting on the Euclidean space $\C^2$. Then $\|A\| = 1$,
$$ p(A)= \left[\begin{array}{c@{\quad}c}
p(a) & ~b \frac{p(a)-p(c)}{a-c} \\
0 & ~p(c)\end{array}\right]
$$
and the von Neumann inequality $\|p(A)\| \le \|p\|_{\overline{\D}}$ implies that
$$ \left| \frac{p(a)-p(c)}{a-c}\right|^2 \le \frac{1 - |p(a)|^2}{1-|a|^2} \frac{1 - |p(c)|^2}{1-|c|^2} .$$
Using the identity $|1-u\overline{v}|^2 = (1-|u|^2)(1-|v|^2) + |u-v|^2$, this can be written as
$$ \left|\frac{p(a)-p(c)}{1-\overline{p(c)}p(a)}\right| \le \left|\frac{a-c}{1-\overline{c}a}\right| .$$
 \end{enumerate}

\vspace*{1pc}
\section{The Multidimensional von Neumann Inequality}\label{sMvN}

\df

\smallskip

\noindent
We say that the {\bf multidimensional von Neumann inequality} holds for a fixed
$n$-tuple of commuting operators $A = (A_1,A_2, \ldots, A_n)$ if
$$\|p(A_1,\ldots, A_n)\| \le \|p\|_{\overline{\D}^n} $$
for every polynomial $p$ in $n$ (commutative) variables.

\smallskip

\fc

\smallskip

We use \cite{NikolskiEasy} as a general reference.

\begin{enumerate}

\item The multidimensional von Neumann inequality holds in the following situations:

\begin{enumerate}
\item \cite{Ando} for a pair of commutative Hilbert space contractions ($n=2$).
\item for a commutative family of isometries.
\item for a family of \emph{doubly commuting} (i.e., $A_iA_j = A_jA_i$ for all $i$ and $j$
and $A_i^*A_j = A_jA_i^*$ whenever $i\neq j$) contractions.
\item \cite{brehmer} for a commutative family $A$
such that $\sum_{j=1}^n \|A_j\|^2 \le 1$.
\end{enumerate}

\item \cite{AglerMcC} Distinguished varieties: let $A_1$ and $A_2$ be two commuting contractive matrices, neither of which has eigenvalues of modulus one. Then there is a polynomial $q\in \C[z,w]$ such that the algebraic set $V= \{(z,w) \in \D^2: q(z,w) = 0\}$ verifies $\|p(A_1,A_2)\| \le \|p\|_V$ for any polynomial $p$ in two variables and
$\overline{V}\cap \partial(\D^2) = \overline{V}\cap (\partial \D)^2$ (the variety exits the bidisk through the
distinguished boundary).

\item Extremal $n$-tuples.

\begin{enumerate}

\item \cite{DruryPAMS,Popescu,arvesonIII} The Drury-Arveson space and the von Neumann inequality of Drury-Popescu-Arveson:
let $\mathbb{B}^n = \{z\in \C^n : \sum_{j=1}^n|z_j|^2 < 1\}$ be the unit open
ball in $\C^n$. Let $DA_n$ be the Drury-Arveson space of all power series $g$ such that
$$g = \sum_{\alpha \ge 0} a_{\alpha}z^{\alpha}, \quad \|g\|_{DA_n}^2 =
\sum_{\alpha \ge 0} |a_{\alpha}|^2\frac{\alpha !}{|\alpha| !} < \infty ,$$
where $\alpha = (\alpha_1, \dots, \alpha_n) \in \N^n$, $|\alpha| = \sum_{j=1}^n \alpha_j$
and $\alpha ! = \prod_{j=1}^n (\alpha_j !)$. The Drury-Arveson space $DA_n$ can be
also regarded as the reproducing kernel Hilbert space with kernel
$$ k_{\lambda}(z) = \frac{1}{1 - \sum_{j=1}^nz_j\overline{\lambda_j}}, \quad z,\lambda \in \mathbb{B}^n. $$
Let $S_jg(z) = z_jg(z)$, where $z = (z_1,\dots,z_n)$ and $g\in DA_n$.
Let $A = (A_1, \dots, A_n)$ be a commutative $n$-tuple on a Hilbert space
$H$ such that $I \succeq \sum_{j=1}^n A_jA_j^{\ast} $, or equivalently,
$$ \left\| \sum_{j=1}^n A_jx_j\right\|^2 \le  \sum_{j=1}^n \left\| x_j\right\|^2$$
for all $x_j \in H$. Then
$$ \|p(A_1,\ldots, A_n)\| \le \|p(S_1,\ldots, S_n)\| $$
for every polynomial $p$.

\item \cite{DruryPAMS} A dual version: 
let $S_j^{\ast}$, $j= 1, \dots, n$, be the
backward shift operators on the Cauchy dual $DA_n^{\ast}$ of the Drury-Arveson space. This means that
$g = \sum_{\alpha \ge 0} a_{\alpha}z^{\alpha}$ with $\|g\|_{DA_n^{\ast}}^2 =
\sum_{\alpha \ge 0} |a_{\alpha}|^2\frac{|\alpha| !}{\alpha !} < \infty ,$ and
$a_{\alpha}(S_j^{\ast}g) = a_{\alpha + \delta_j}(g)$ for all $\alpha \ge 0$,
where $\delta_j = (\delta_{jk})_{1\le k \le n}$.
Let $A = (A_1, \dots, A_n)$ be a commutative $n$-tuple on a Hilbert space
$H$ such that $I \succeq \sum_{j=1}^n A_j^{\ast}A_j$, or equivalently,
$$ \sum_{j=1}^n \left\|A_jx\right\|^2 \le \left\| x\right\|^2$$
for all $x \in H$. Then
$$ \|p(A_1,\ldots, A_n)\| \le \|p(S_1^{\ast},\ldots, S_n^{\ast})\| $$
for every polynomial $p$.
\end{enumerate}
\item \cite{Mlak1971} Mlak's von Neumann inequality with operator coefficients:
let $C_k : H \mapsto H$ be bounded linear operators, $k=0,\ldots,n$, and $A \in \LH$ be a
contraction which double-commutes with the $C_k$'s, i.e. $AC_k = C_kA$ and $AC_k^\ast = C_k^\ast A$ for every $k$.
Then
$$ \left\| \sum_{j=0}^nC_jA^j\right\| \le \sup_{|z|\le 1}\left\| \sum_{j=0}^n C_jz^j\right\| .$$

\item \cite{Lubin} Let $(A_1, \ldots, A_n)$ be commuting contractions on a Hilbert space $H$.
Then the polydisk $\{(z_1, \ldots , z_n) : |z_i| < \sqrt{n}, 1\le j \le n\}$ is a spectral
set for $(A_1, \ldots, A_n)$.

\item \cite{Hirschfeld} The \emph{Poisson radius}:
let $A = (A_1,\ldots, A_n)$ be an $n$-tuple of commuting operators
on a Hilbert space $H$, with spectra included in $\overline{\D}$. Define
$$P_j(rA_j, \zeta_j) = \mathrm{re}\left( \left( (\zeta_j I + rA_j)(\zeta_j I - rA_j)^{-1} \right) \right)$$
for $0 \le r < 1$, $\zeta_j \in \T$, $1\le j \le n$, and the Poisson radius $r_P(A)$ of $A$ as the supremum of
$r \in [0,1)$ such that
$$\frac{1}{n!}\sum_{\sigma} P_{\sigma(1)}(rA_{\sigma(1)}, \zeta_{\sigma(1)})
\ldots P_{\sigma(n)}(rA_{\sigma(n)}, \zeta_{\sigma(n)})$$
is a positive operator for every $\zeta\in \T^n$, where $\sigma$ runs over all permutations of $(1,\ldots,n)$.
Then $0 < r_P(A) \le 1$ and $\|f(r_P(A)A)\| \le \|f\|_{\overline{\D}^n}$ for every polynomial $f$.

\item \cite{AglerVol}  The Schur-Agler class: let $f$ be an analytic function of $n$
complex variables $\lambda = (\lambda^1,\cdots, \lambda^n)$.  Then $f(rA_1,\dots, rA_n)$ has norm at most $1$ for any $r<1$ and any collection
of $n$ commuting contractions $(A_1, \dots, A_n)$ on a Hilbert space if and only if there are auxiliary Hilbert spaces $H_j$, $1\le j\le n$, and an isometry $V \in \mathcal{L}( \C\oplus H_1 \oplus \cdots \oplus H_n)$ such that, if $H =   H_1 \oplus \cdots \oplus H_n$, $V$ is written with respect to $\C\oplus H$ as $V= \left[\begin{array}{c@{\quad}c}
A & B \\
C & D \end{array}\right] ,$
and $\mathcal{E}_{\lambda} = \lambda^1I_{H_1}\oplus \cdots \lambda^n I_{H_n}$, then $f(\lambda) = A+B\mathcal{E}_{\lambda}(I_H-D\mathcal{E}_{\lambda})^{-1}C$.

\item \cite{bozejko} Bo\.zejko's von Neumann inequality for non-commuting tuples:
let $A_k$, $1\le k \le n$, be (not necessarily commuting)
contractions on $H$. Let $f = f(x_1, \ldots, x_n)$ be a polynomial in the noncommutative
indeterminates $x_1, \ldots, x_n$. Then
$$\|f(A_1, \ldots, A_n)\| \le \sup \{\|f(U_1, \ldots, U_n)\| :
U_j \textrm{ unitary matrices on } \C^{m\times m}, m \in \N \} .$$
\end{enumerate}

\smallskip

\noindent{\bf Open Problems:}
  \label{c107sMvNopen}


\begin{enumerate}
\item  \cite{DixonJLMS} It is not known if for each $n$ there exists a finite constant $C_n$ such that
for any commuting contractions $A_1, \dots, A_n$ and any polynomial $f$ in $n$ variables one has
$$ \|f(A_1,\ldots, A_n)\| \le C_n\|f\|_{\overline{\D}^n} .$$
It is generally believed that such a constant $C_n$ does not exist.
One knows that $C_n$ must increase faster than any power of $n$.
\end{enumerate}

\medskip
\ex

\begin{enumerate}
\item \cite{Varopoulos} The multidimensional von
Neumann inequality
fails in general
for $n\ge 3$ and matrices $A_1,\ldots, A_n\in \C^{d\times d}$.
\item  \cite{Varopoulos} The multidimensional von Neumann inequality can
fail with $n=3$ and $d=5$. The three matrices $A_1,A_2,A_3\in \C^{5\times 5}$ are commuting contractions (with respect to the Euclidean norm), and the
polynomial $p(z_1,z_2,z_3) = z_1^2+ z_2^2 + z_3^2 -2z_1z_2 -2z_1z_3 -2z_2z_3$
satisfies $\|p\|_{\T^3} = 5$ and $\|f(A_1, A_2, A_3)\| > 5$:
$$ A_1 =
\left[\begin{array}{c@{\quad}c@{\quad}c@{\quad}c@{\quad}c}
0& 0 & 0 & 0 & 0\\
1& 0 & 0 & 0 & 0\\
0& 0 & 0 & 0 & 0\\
0& 0 & 0 & 0 & 0\\
0& \frac{1}{\sqrt{3}} & -\frac{1}{\sqrt{3}} & -\frac{1}{\sqrt{3}}~ & 0
\end{array}\right],
~
A_2 =
\left[\begin{array}{c@{\quad}c@{\quad}c@{\quad}c@{\quad}c}
0& 0 & 0 & 0 & 0\\
0& 0 & 0 & 0 & 0\\
1& 0 & 0 & 0 & 0\\
0& 0 & 0 & 0 & 0\\
0~& -\frac{1}{\sqrt{3}}~ & \frac{1}{\sqrt{3}} & -\frac{1}{\sqrt{3}}~ & 0\end{array}\right]
,
$$

$$
A_3 =
\left[\begin{array}{c@{\quad}c@{\quad}c@{\quad}c@{\quad}c}
0& 0 & 0 & 0 & 0\\
0& 0 & 0 & 0 & 0\\
0& 0 & 0 & 0 & 0\\
1& 0 & 0 & 0 & 0\\
0& -\frac{1}{\sqrt{3}} & -\frac{1}{\sqrt{3}}~ & \frac{1}{\sqrt{3}}~ & ~0
\end{array}\right]
.$$
\item \cite{CrabbDavie1975} Denote by $\mathbf{e}_j$, $1\le j \le 8$, the vectors of the standard orthonormal
basis of $\C^8$. Let $p(z_1,z_2,z_3) = z_1z_2z_3 -z_1^3-z_2^3-z_3^3$. There exist
three commuting matrices $A_k\in \mathbb C^{8 \times 8}$, $k=1,2,3$, such that $\|A_i\| \le 1$ and
$ \|p(A_1,A_2,A_3)\| \ge 4 > \|p\|_{\T^3}$.
The contractions $A_i$ are acting on the orthonormal basis as follows:\vspace{-5pt}
\renewcommand{\arraystretch}{1.2}
$$
 \begin{array}{ccccccccccccc}
A_1 : \mathbf{e}_1 & \longrightarrow \mathbf{e}_2 & \longrightarrow (-\mathbf{e}_5) & \longrightarrow (-\mathbf{e}_8) & \longrightarrow 0, & \,
& \mathbf{e}_3 & \longrightarrow \mathbf{e}_7 & \longrightarrow 0, &
\, & \mathbf{e}_4 & \longrightarrow \mathbf{e}_6 & \longrightarrow 0 \\
A_2 : \mathbf{e}_1 & \longrightarrow \mathbf{e}_3 & \longrightarrow (-\mathbf{e}_6) & \longrightarrow (-\mathbf{e}_8) & \longrightarrow 0, &
\, & \mathbf{e}_2 & \longrightarrow \mathbf{e}_7 & \longrightarrow 0, &
\, & \mathbf{e}_4 & \longrightarrow \mathbf{e}_5 & \longrightarrow 0 \\
A_3 : \mathbf{e}_1 & \longrightarrow \mathbf{e}_4 & \longrightarrow (-\mathbf{e}_7) & \longrightarrow (-\mathbf{e}_8) & \longrightarrow 0, &
\, & \mathbf{e}_2 & \longrightarrow \mathbf{e}_6 & \longrightarrow 0, &
\, & \mathbf{e}_3 & \longrightarrow \mathbf{e}_5 & \longrightarrow 0
\end{array}
$$
\renewcommand{\arraystretch}{1}

\item \cite{Varopoulos} Given $K> 0$, there exist a positive integer $n$,
commuting operators $A_1,\ldots, A_n$ and a polynomial $p$
such that $\sum_{j=1}^n \|A_j\|^2 \le 1$ and\vspace{-5pt}
$$ \|p(A_1,\ldots, A_n)\| > K \sup \left\{ |p(z_1,\ldots,z_n)| : \sum_{j=1}^n |z_j|^2 \le 1 \right\} .$$
\end{enumerate}

\vspace*{1pc}
\section{Dilations, Complete Bounds and Similarity Problems}\label{sDilat}
\smallskip

\df

\smallskip

\noindent
For a given closed set $X$ of the complex plane,
we say that $A \in \LH$ has a {\bf normal} $\partial X$-{\bf dilation} if there exist a Hilbert
space $\mathcal H$
containing $H$ and a normal operator $N$ on $\mathcal H$ with $\sigma(N)\subset \partial X$ so that
$$ f(A) = P_Hf(N)\mid_H$$ for every rational function $f$ with
poles off $X$. Here $P_H$ is the orthogonal projection of $\mathcal{H}$ onto $H$.

If $X = \overline{\D}$, then $N$ is a unitary operator and we say that $A$ has a {\bf unitary (strong) dilation}.
Notice that in this chapter, contrary to Chapter 18, a unitary (strong) dilation $A$ is a
common dilation of all powers of $A$.

For a fixed $s > 0$, we say that
$A \in \LH$ has a $s$-{\bf unitary dilation} if there exists a Hilbert space $\mathcal H$
containing $H$ and a unitary operator $U$ on $\mathcal H$ such that
$$ A^n = s P_H U^n\mid_H , \quad (n\ge 1).$$
The last two definitions agree in the case $s = 1$.

We denote by $M_n(\functionsR (X))$ the algebra of $n$ by $n$ matrices with entries from $\functionsR (X)$.

Considering the (spectral) H\"older $2$-norm  for matrices in $\mathbb C^{n \times n}$, we can endow
$M_n(\functionsR (X))$ with the norm
$$\|\left(f_{ij}\right)_{1\le i,j \le n}\|_X = \sup \{ \|\left(f_{ij}(x)\right)_{1\le i,j \le n}\| :
x \in X\} = \sup \{ \|\left(f_{ij}(x)\right)_{1\le i,j \le n}\| : x \in \partial
X\}. $$ In a similar fashion we endow $M_n(\mathcal{L}(H))$ with
the norm it inherits by regarding an element $(A_{ij})_{1\le i,j \le n}$ in
$M_n(\mathcal{L}(H))$ as an operator acting on the direct sum of
$n$ copies of $H$.

For a fixed constant $K > 0$, the set $X$ is
said to be a {\bf complete} $K$-{\bf spectral} set for $A$ if
$\sigma(A) \subset X$ and the inequality $\|(f_{ij}(A))_{1\le i,j \le n}\| \le
K\|(f_{ij})_{1\le i,j \le n}\|_X$ holds for every matrix $(f_{ij}) \in M_n(\functionsR
(X))$ and every $n$.

A {\bf complete spectral} set is a
complete $K$-spectral set with $K=1$.

We also say that
$A$ is
{\bf power bounded} if $\sup_n\|A^n\| < \infty$.

Two Hilbert space operators $A$ and $B$ are
said to be {\bf similar} if there exists an invertible operator $L$ such that $B = L^{-1}AL$.

\bigskip

\fc

\smallskip

\noindent
All the following facts except those with a specific reference can be found
in
\cite{FoiasSzNagy, PisierLNM, Pa02}. We denote by $A$ a Hilbert space
operator and by $X$ a closed subset of $\C$.

\begin{enumerate}

\item \cite{arvesonI,arvesonII, paulsenPAMS} $A$ has a normal $\partial X$-dilation if and
only if $X$ is a complete spectral set for $A$.

\item $X$ is completely $K$-spectral for $A$ if and only if $X$ is completely spectral
for an operator $B\in \LH$ similar to $A$, say, $B = L^{-1} A L$, with $\| L^{-1}\| \, \|L\| \leq  K$.

\item Each completely spectral set for $A$ is spectral. Conversely,  a spectral set for $A$ is
completely spectral in the following situations:

\begin{enumerate}
\item  if $X$ is a closed disk;

\item \cite{agler} if $X$ is an annulus;

\item \cite[Theorem 4.4]{Pa02} if $\C\setminus X$ has only finitely many components and
the interior of $X$ is simply connected or, more generally, if $\functionsR(X) + \overline{\functionsR(X)}$ is dense in $\functionsC (\partial X)$.

\item \cite{aglerMem,DritschelMcC}
There is a closed set $X$ in $\C$ having ``two holes'' and an operator $A$
such that $X$ is spectral for $A$ but not completely spectral.
\cite{Pickering} More generally, such a counterexample exists whenever $X$ is a symmetric domain in $\C$ with $n$ holes, $2\le n < \infty$.
\item \label{c107.s4.fa2}\cite{dopa,Pa02} 
     Let the boundary of the compact $X \subset \mathbb C$ consist of $n+1$ disjoint Jordan curves. If $X$ is a spectral set then it is a complete $(2n+1)$-spectral set.

\end{enumerate}

\item
Each completely $K$-spectral set for $A\in \LH$ is $K$-spectral.
\cite{PauJOT} The converse is true whenever $A$ is a $2\times 2$ matrix.
\cite{PisierJAMS} However, there exists an example where $\overline{\mathbb D}$ is $K$-spectral for $A\in \LH$ but not completely $K'$-spectral for any $K'$.

\item A combination of these facts yields several corollaries which historically came first:

\begin{enumerate}
\item (Sz.-Nagy dilation theorem) Every Hilbert space contraction has a unitary dilation.

\item (Berger dilation theorem) Every numerical radius contraction has a unitary $2$-dilation.

\item (Paulsen criterion) An operator $A$ is similar to a contraction if and only if it has the closed unit disk as a complete $K$-spectral set.
\end{enumerate}

\item Operators of class C$_{s}$:

\begin{enumerate}

\item \cite{FoiasNagySzeged,OkuboAndo} Let $s > 0$. Every operator of class
C$_{s}$ has a unitary $s$-dilation.

\item Let $s > 0$. If $A$ is of class C$_{s}$, then
$\|f(A)\| \le  \max \{ |s \cdot f(z)+(1-s) \cdot f(0)|  : |z| \le 1\}
$ for every polynomial $f$.

\item Let $H$ be a complex Hilbert space of dimension $\ge 2$. Then the class C$_{s}$ increases with $s$:
we have  C$_{s} \subset $  C$_{s'}$ and C$_{s} \neq $  C$_{s'}$ for $0 < s < s'$.
The set of operators acting on $H$ which belong to one of the classes C$_{s}$, for some $s > 0$, is
dense in the strong operator topology in the set of all power bounded operators.

\end{enumerate}

\item Unitary dilations for $n$-tuples:

\begin{enumerate}

\item \cite{Ando} Every pair of commuting contractions on a Hilbert space has a pair of commuting unitary dilations.

\item \cite{GasparRacz} Every
$n$-tuple $A = (A_1,\ldots,A_n) \in \LH^n$ which is {\em cyclic commutative}, i.e.
$$A_1A_2\ldots A_n = A_nA_1\ldots A_{n-1} = \ldots = A_2A_3 \ldots A_nA_1 ,$$
has a cyclic commutative dilation to an $n$-tuple of unitaries.

\item \cite{Opela}
Let $G$ be an {\em acyclic} graph on $n$ vertices $\{1, 2, \dots , n\}$
(this means that it does not contain a cycle as a subgraph).
Let $A=(A_1,A_2, \dots,A_n)$ be an $n$-tuple of contractions on a Hilbert space that {\em commute according
to} $G$, that is $ A_iA_j = A_jA_i$ whenever $(i,j)$ is an edge of $G$.
Then there exists an $n$-tuple $U$ of
unitaries on a larger Hilbert space that commute according
to $G$ such that $U$ dilates $A$. This property may fail if $G$ contains a cycle.
\end{enumerate}

\item Power boundedness and
similarity to a contraction (i.e., $\overline{\mathbb D}$ is completely $K$-spectral):

\begin{enumerate}

\item 
\cite{Rota,Herrero,Voiculescu} Rota type theorems: If $A\in \LH$ with $\sigma(A) \subset \D$,
then $A$ is similar to a contraction. More general results are true for some closed sets $X\subset \C$ and
for operators $A$ with $\sigma(A) \subset X$.

\item \cite{BadeaJOT} A Banach space Rota theorem and Matsaev inequality:
Let $E_1$ be a Banach space and
suppose that $A \in \mathcal{L}(E_1)$ has $\sigma(A) \subset \D$. Then, for every $p > 1$, there exists
	a Banach space $E_2$ which is a quotient of the set $\ell^{p}(E_1)$ of $E_1$-valued sequences
and an
	isomorphism $L : E_2 \to E_1$ such that, if $B = L^{-1}AL \in
	\mathcal{L}(E_2)$, then
$$
	\|f(B)\|_{\mathcal{L}(E_2)} \leq \|f(S)\|_{\mathcal{L}(\ell^{p}(E_1))}
$$
	for each polynomial $f$ and, even more generally,
	$$
       \|(f_{ij}(B))_{1\le i,j\le n}\|_{\mathcal{L}(E_2^n)}
	     \leq \|(f_{ij}(S))_{1\le i,j\le n}\|_{\mathcal{L}(\ell^{p}(E_1)^n)}
    $$
	for all matrices of polynomials.

\item 
If $A$ is compact and power bounded, then $A$ is similar to a contraction.

\item \cite{Foguel,Lebow, PisierJAMS} There is a Hilbert space operator $A$ which is power bounded but for which the closed
    unit disk is not $K$-spectral for any $K$, and thus $A$ is not similar to a contraction. There is a
Hilbert space operator $A$ which is not similar to a contraction but for which the closed unit disk is $K$-spectral
for $A$, for some $K$.

\item \cite{Peller82,BourgainIsrael} Let $A\in \LH$ with $M = \sup_n\|A^n\| < \infty$ and
let $p$ be a polynomial of degree $d \ge 2$. Then $\|p(A)\| \le M^2\, (\log d) \, \|\, p \, \|_{\overline{\D}}$. A similar
result holds for norms of matrices of polynomials, and the $\log d$ term in the inequality is the
best one may hope for.

\item \cite{BourgainIsrael} Bourgain's estimate for matrices similar to contractions: If $A$ is a matrix
such that, with the (spectral) H\"older $2$-norm,
$\|f(A)\| \le C\|f\|_{\overline{\D}}$ for any polynomial $f$, then there is an invertible matrix $
L$ such that $\|L^{-1}AL\| \le 1$ and $\|L^{-1}\| \, \|L\| \le KC^4 \log (n+1)$, where $K$ is
a numerical constant independent of $n$.

\item \cite{BercoviciPrunaru} Suppose that $\overline{\mathbb D}$ is
 $K$-spectral for $A\in \LH$. Then
there exist Hilbert spaces $H_1,H_2$, contractions $A_1 \in \mathcal{L}(H_1)$,
$A_2\in \mathcal{L}(H_2)$ and injective linear operators $X_1 : H_1 \mapsto H, X_2 :H \mapsto H_2$ with dense ranges such
that $X_1A_1 = AX_1$ and $A_2X_2 = X_2A$.
See also Example~\ref{sDilat}.\ref{exMullerTomilov}.
\end{enumerate}
\end{enumerate}

\smallskip

\noindent{\bf Open Problems:}
  \label{c107sDilatopen}

\begin{enumerate}
\item \cite{FoiasSzNagy} What is the obstruction to an
$n$-tuple of commuting Hilbert space contractions having
commuting unitary dilation?
\end{enumerate}

\medskip


\ex
\smallskip
\begin{enumerate}

\item 
The (unilateral) backward shift operator $S^{\ast}\in \mathcal L(\ell^2)$ defined by $(S^{\ast}x)_j = x_{j+1}$, $j\in \mathbb N$, is easily seen to be a contraction, with spectrum $\sigma(S^{\ast})=\overline{\mathbb D}$. Hence $X=\overline{\mathbb D}$ is a (completely) spectral set for $S^{\ast}$. For quadratic summable sequences $x$ indexed by $\mathbb Z$, one defines the bilateral backward shift operator $B$ through $(Bx)_j = x_{j+1}$, $j\in \mathbb Z$, which is clearly unitary. Since $(B^nx)_j = x_{j+n}$, $j\in \mathbb Z$, $n\geq 1$, we see that $B$ is a $1$-unitary dilation and thus a unitary dilation of $S^{\ast}$.

\item The notion of (weak or strong) unitary dilation is a nice illustration of why it is
important in linear algebra to have sometimes recourse to infinite dimensions:
we say that $A \in \mathbb C^{n\times n}$ is imbedded in $B\in\mathbb C^{m \times m}$ if
$B = \left[\begin{array}{c@{\quad}c} A  &  *  \\  *  &  * \end{array}\right]
$. Given a
contraction $A\in \mathbb C^{n\times n}$, is there a unitary $B\in \mathbb C^{m\times m}$
such that $p(A)$ is imbedded in $p(B)$ for any polynomial $p$ of degree at most $k$?
This property, related with the exactness property for Krylov spaces, reduces to
(strong) unitary dilations for $k=\infty$ considered above, and
for $k=1$ to the (weak) unitary dilation considered in
Section 18.6.
 The latter problem has the solution \cite{halmosBR}
    $$
       B=\left[\begin{array}{c@{\quad}c}A  &  (I-A A^*)^{1/2}  \\ (I-A^*A)^{1/2}  &  - A^* \end{array}\right].
    $$
    Eg\'ervary \cite{Egervary} showed that such an imbedding is always possible with $m=(k+1)n$:
for example, one may imbed the $n$-dimensional shift in an unitary circulant of order
$m=(k+1)n$. It is also known that in general such an imbedding for finite $m$ is
impossible for $k=\infty$, see also \cite{LeSh,McSh} for
more modern aspects of this question.

\item \cite{Parott} 
Let $U$ and $V$ be two contractions in $\LH$ such that $U$ is unitary
and $UV \neq VU$. We define three commuting contractions in $\mathcal{L}(H\oplus H)$ by defining\vspace{-5pt}
$$A_1=\left[\begin{array}{c@{\quad}c}0 & 0 \\I & 0 \end{array}\right], 
A_2=\left[\begin{array}{c@{\quad}c}0 & 0 \\U & 0 \end{array}\right], A_3=\left[\begin{array}{c@{\quad}c}0 & 0 \\V & 0 \end{array}\right].$$
Then the commuting triplet $A= (A_1,A_2,A_3)$ verify the multidimensional
von Neumann inequality but do not possess a commuting triplet of (strong) unitary dilations.

\item \label{diagonalizable}
Let $A \in \mathbb C^{n\times n}$ be diagonalizable. Then $A$ is power bounded if
and only if $\sigma(A)\subset \overline{\mathbb D}$, and in this case
    $A$ is similar to the contraction given by the Jordan canonical form.
    Here $\overline{\mathbb D}$ is (completely) $K$-spectral, with $K$ the condition number of the matrix of eigenvectors.

    If $A \in \mathbb C^{n\times n}$ is not diagonalizable,
then $\sigma(A)\subset \mathbb D$ still implies
(through, e.g., the Cauchy integral formula and the pseudo-spectrum)
that $\overline{\mathbb D}$ is completely $K$-spectral for some $K$.
However, in general, the Jordan canonical form is no longer a contraction.
\item \label{exMullerTomilov} \cite{MuTo} 
Let $T \in \LH$ and $S \in \mathcal{L}(K)$ be Hilbert space operators.
We say that $T$ is a \emph{quasiaffine transform} of $S$ if there
exists an injective operator $A:H \mapsto K$ with dense range such that $AT = SA$.
We say that $T$ is \emph{quasisimilar} to $S$ if each operator is a quasiaffine transform of the other.

There exists a power bounded operator
on a Hilbert space which is not quasisimilar to a contraction.

\item \cite{PisierJAMS, DavidsonPaulsen, BadeaJOT, BadeaPaulsen, Ricard}
Let $\alpha=(\alpha_0,\alpha_1,\dots)$ be a sequence in $\ell^2$ and
set
$$
R(\alpha) =
   \left[\begin{array}{c@{\quad}c}
S^{\ast}& Y(\alpha)\\0& S\end{array}\right]
\in \mathcal{L}(\ell^2(H)\oplus \ell^2(H)),$$
where $S$ is the shift on the Hilbert space $\ell^2(H)$ of $H$-valued square summable
sequences, $H$ is of infinite dimension, and
$
 Y(\alpha) = \Bigl[ \alpha_{i+j}C_{i+j} \Bigr]_{i,j\ge 0},
$
where the
$C_j$'s are operators verifying the \emph{canonical anticommutation relations}
$ C_iC_j + C_jC_i = 0$ and $C_iC_j^{\ast} + C_j^{\ast}C_i = \delta_{ij} I .$

Then the operator $R({\alpha})$ is polynomially bounded if and only if
$\sup_{k\ge0} (k+1)^2 \sum_{i\ge k}|\alpha_i|^2$
is finite, and $R({\alpha})$ is similar
to a
contraction if and only if $\sum_{k\ge0} (k+1)^2 |\alpha_k|^2$ is finite.


\end{enumerate}

\vspace*{1pc}
\section[Intersections of Spectral and $K$-Spectral Sets]{Intersections of Spectral and {\large $K$}-Spectral Sets
}\label{s4}

In this part we will discuss intersections of spectral sets, including the annulus problem of Shields \cite{shields}. It is known that the intersection of two spectral sets is not necessarily a spectral set, see \cite{WilliamsSzeged, misra, Pa02} and Example~\ref{s4}.\ref{c107.s4.ex1} below. However, the same question for $K$-spectral sets remains open, though the above examples indicate that one may not use the same constant.
We refer to \cite{paulsenK}
and the book \cite{Pa02} for modern surveys of known properties of $K$-spectral and complete
$K$-spectral sets.

       \begin{figure}[ht!]
\centerline{\scalebox{.75}{\includegraphics{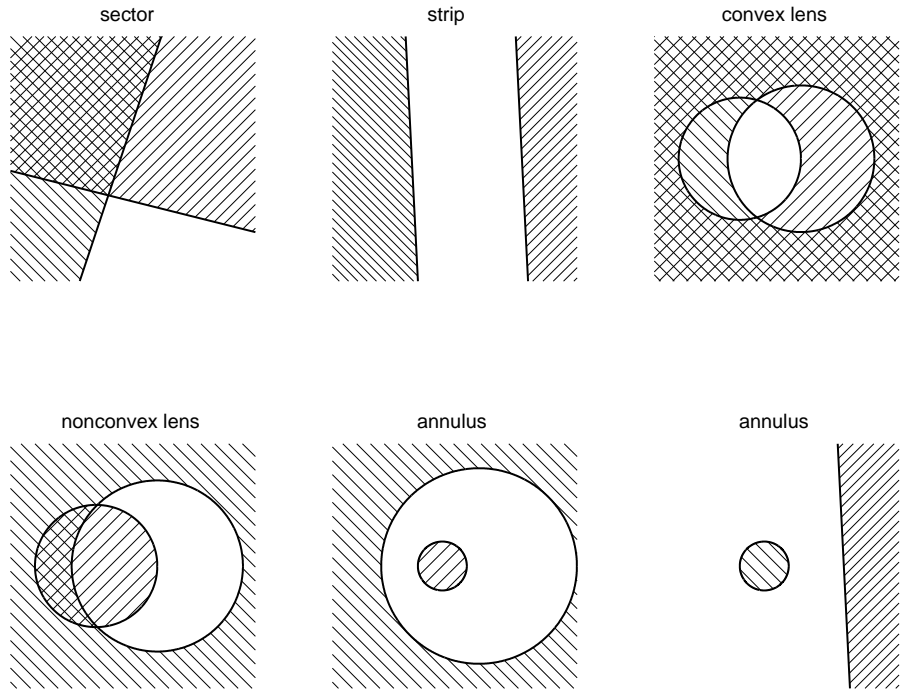}}}
\caption{\label{c107fig1}Six different configurations (in white) of intersections of two disks of the Riemann sphere.}
\end{figure}

\fc


\begin{enumerate}

\item \cite{dopa,Pa02}
     If two $K$-spectral sets have disjoint boundaries, then their intersection is a $K'$-spectral set for some $K'$.

\item \cite{stampfliI} The intersection of a simply connected spectral set $X$ of $A \in \LH$ whose interior has finitely many components with the closure a simply connected open set $G$ containing the spectrum of $A$ is a $K'$-spectral set for $A$ for some $K'$.

    Weaker versions of this statement concerning the connectivity of $X$ and/or $G$ have been given in \cite{stampfliII,lewis}.
\item \cite{lewis} The
      intersection of a (complete) $K$-spectral set for a bounded linear operator $A$ with the closure of any open set containing
      the spectrum of $A$ is a (complete) $K'$-spectral set for $A$ for some $K'$.
\item \cite{BBC09} Let $A\in \mathcal L(H)$, and consider the
      intersection $X=D_1\cap D_2\cap\dots\cap D_n$ of $n$ disks of the Riemann sphere $\overline \C$, each of them being spectral for $A$. Then $X$ is a complete $K$-spectral set for $A$, with a constant $K\leq n+n(n\!-\!1)/\sqrt3$.
      \begin{enumerate}
        \item \cite{Cr07} If in this result we add the requirement that the disks $D_j$ and thus $X$ are convex, then
            $X$ is a complete $11.08$-spectral set for $A$.
 
        \item For $n=2$ we obtain the constant $K=2+2/\sqrt{3}$ for various configurations as shown in Figure~\ref{c107fig1}, in particular for a strip/sector obtained by the intersection of two half-planes and discussed in Example~\ref{s5}.\ref{c107s5fa3},
            or the lens-shaped intersection of two disks \cite{BC06}.
      \end{enumerate}
    \item \cite{shields} For $R > 1$, consider the annulus $X = X(R)= \{z \in \C : R^{-1} \le |z| \le R\}$, and denote by $K(R)$ (and by $K_{cb}(R)\geq K(R)$, respectively), the smallest constant $K$ such that $X$ is a $K$-spectral set (and a complete $K$-spectral set,
        respectively) for any invertible $A \in \mathcal{L}(H)$ verifying
        $\|A\| \le R$ and $\|A^{-1}\| \le R$. Then
        \begin{enumerate}
        \item \label{c107.s4.fa6a} \cite{shields} $K(R) \leq 2 + \sqrt{\frac{R^2+1}{R^2-1}}$;
        \item \label{c107.s4.fa6b}  \cite{BBC09} $K_{cb}(R) \leq 2 + \frac{R+1}{\sqrt{R^2+R+1}}\leq 2 + \frac{2}{\sqrt{3}}$. The first upper bound is sharper than \ref{s4}.\ref{c107.s4.fa6a} for $R\leq 3.1528$;
        \item \label{c107.s4.fa6c} \cite{BBC09} $K(R)\geq \frac43$.
        \item \label{c107.s4.fa6d} \cite{BBC09}
        $K_{cb}(R)\leq
\max \bigl\{ 3 , 2+\sum_{n=1}^\infty \frac{4}{R^{2n} + 1} \bigr\}$,
being sharper than \ref{s4}.\ref{c107.s4.fa6a} for $R \geq 1.8544$, and sharper than \ref{s4}.\ref{c107.s4.fa6b} whenever $R\geq 1.9879$.
        It follows in particular that $K(R) \leq 3$ for $R\geq 2.0953$.
\item \label{c107.s4.fa6e} \cite{Cr11}
        $K_{cb}(R)\leq 2 + \frac{1}{\pi} \int_0^\pi \Bigl| \frac{R^2+\exp(i\theta)}{R^2-\exp(i\theta)} \Bigr| \, d\theta$, being
        always sharper than \ref{s4}.\ref{c107.s4.fa6a}, sharper than \ref{s4}.\ref{c107.s4.fa6b} for $R\geq 1.6405$, and sharper than \ref{s4}.\ref{c107.s4.fa6d} whenever $R\leq 2.0462$.
%
        \end{enumerate}
\end{enumerate}
\medskip
\ex

\smallskip

\noindent
  Consider the matrix $$
           A=\left[\begin{array}{r@{\quad}r}
1 &\gamma\\0 &1\end{array}\right]\in \mathbb C^{2\times 2}, \quad \gamma=R-\frac{1}{R}, \quad R>1,
     $$
     with $\| A \| = \| A^{-1} \| = R$ such that both sets $\{|z|\leq R\}$ and $\{|z|\geq 1/R\}$ are spectral for $A$, and consider the intersection $X(R)= \{z \in \C : R^{-1} \le |z| \le R\}$.
 \begin{enumerate}
 \item \label{c107.s4.ex1}
     \cite{misra} The example of the function $f(z)=z-1/z$ which verifies $\| f(A) \|/\| f \|_{X(R)}=2\frac{R^2-1}{R^2+1}$ shows that $X(R)$ is not a spectral set for $R > \sqrt{3}$.
 \item We get a sharper statement for the function $f(z)=g(z)-g(1/z)$, $g(z)=R \frac{z-1}{R^2-z}$, leading to
     $$
         \| f(A) \| = 2 , \quad \| f \|_{X(R)} = \frac{1+R^2+2R}{1+R^2+R}< \frac{4}{3} .
     $$
     Thus $X(R)$ is even not $\frac{3}{2}$-spectral for $A$ for any $R>1$. Compared with Fact~\ref{s4}.\ref{c107.s4.fa6c} we thus have shown the improved lower bound $K(R)\geq 2 \frac{1+R^2+R}{1+R^2+2R}>\frac{3}{2}$.
 \item
 Let the boundary of the compact $X \subset \mathbb C$ consists of $n+1$ disjoint analytic Jordan curves.
     Then
     $X$ is $K$-spectral for the above matrix $A$ if and only if \cite{misra}
     $$
           K \geq \gamma \, \Gamma(X) , \quad
           \Gamma(X)= \sup \Bigl\{ \frac{|f'(1)|}{\| f \|_{X}} : \mbox{$f$ analytic on the interior of $X$, $f(1)=0$} \Bigr\} .
     $$
     For the annulus $X(R)=\{ z \in \mathbb C : R^{-1}\leq |z|\leq R \}$ of the preceding example, the quantity $\Gamma(X(R))$ is computable \cite{BBC09}, leading to the lower bound of Fact~\ref{s4}.\ref{c107.s4.fa6c}.

     Moreover, given $r\in (1,R)$,  using the Schwarz lemma we also conclude that both sets $\{z: |z|\leq r\}$ and $\{z: |z|\geq 1/r\}$ are $K$-spectral for $A$ with
     $K=(R-R^{-1})/(r-r^{-1})$ but not their intersection $X(r)$, at least for those $r$ not too far from $R$.

\end{enumerate}

\vspace*{1pc}
\section[The Numerical Range as a $K$-Spectral Set]{The Numerical Range as a {\large $K$}-Spectral Set}\label{s5}
Recall that the numerical range is given by $W(A) :=  \{\langle Ax,x\rangle : \|x\| = 1\}$.
It has been conjectured \cite{FoiasWilliams} that the closure of the numerical range
is a complete $K$-spectral
set for $A$. This was proved in \cite{DelyonDelyon}. Here we report about recent results along these lines, including the
existence of a universal such constant $K$ shown by Crouzeix \cite{Cr07}. For applications in numerical linear algebra of results of this type see for
instance the discussions in \cite{eiermann,greenbaum,tohtrefethen}.

\bigskip

\fc


 \begin{enumerate}
 \item \label{c17.s5.fa1} Estimates depending on the shape of $X$: For every bounded linear operator $A\in \LH$,
every compact convex set $X$ containing the numerical range is a complete $K$-spectral set for $A$, where
     \begin{enumerate}
        \item \cite{DelyonDelyon} $K = 3 + ( 2 \pi \, \mbox{\rm diameter}(X)^2/\,\mbox{\rm area}(X))^3$.
        \item \label{c17.s5.fa1b} \cite{PutinarSandberg}
           $K = 1+2/(1-q(X))$, with $q(X)\in [0,1]$ being C.\ Neumann's configuration constant of $X$ (the oscillation norm of the underlying Neumann-Poincar\'e singular integral operator \cite{kral}). We have $q(X)=1$ if and only if $X$ is a triangle or a quadrilateral.
        \item \label{c17.s5.fa1c} \cite[Theorem~2.3]{BCD06}
           $K = 2+\pi+TV(\log(r))$, if the boundary of $X$ is parametrized by $[0,2\pi] \ni t \mapsto \omega+r(t)e^{it}$, $r(t)\geq 0$.
    \end{enumerate}
  \item \label{c17.s5.fa2} Universal estimates:
    \begin{enumerate}
        \item \cite[Theorem~1]{Cr07} There exists a universal constant $K=K_{Crouzeix}\in [2,11.08]$ such that for every bounded linear operator $A\in \LH$, the closure of the numerical range $\overline{W(A)}$ is a complete $K$-spectral set for $A$.
        \item  \cite[Theorem~1.1]{Cr04a} $W(A)$ is completely $2$-spectral for every $2 \times 2$ matrix $A$.
    \end{enumerate}
  \item Mapping theorems for the numerical range: let $A\in \LH$, and let
      $f$ be analytic on $\mathbb D$, and continuous up to $\partial \D$.
      \begin{enumerate}
        \item \cite{BergerStampfli} $w(A)\leq 1$ if and only if for all $z\in \mathbb C$\vspace{-3pt}
            $$\| A - zI \| \leq 1 + \sqrt{1+|z|^2} = \left\|
            \left[\begin{array}{c@{\quad}c}0 & 2 \\ 0 & 0\end{array}\right] - zI \right\|.$$
        \item \cite{BergerStampfli} If $w(A)\leq 1$, $\| f \|_{\mathbb D}\leq 1$, and $f(0)=0$ then $w(f(A))\leq 1$.
        \item \cite{Drury2008} If $w(A)\leq 1$, $\| f \|_{\mathbb D}\leq 1$, and $|f(0)|<1$, then $W(A)$ is a subset of the convex hull of the union of the disks $\overline{\mathbb D}$ and $\{ z \in \mathbb C: | z-f(0) |\leq 1-|f(0)|^2 \}$. Furthermore, $\|f(A)\| \le \nu(|f(0)|)$, where\vspace{-5pt}
         $$\nu(\alpha) = \left(2-3\alpha^2 +2\alpha^4 +2(1-\alpha^2)(1-\alpha^2+\alpha^4)^{1/2}\right)^{1/2}.$$
        \item \cite{Drury2008} If $w(A) \le 1$, then $w(f(A)) \le \frac{5}{4}\|f\|_{\overline{\D}}$, and $\frac{5}{4}$ is the best possible constant.
        \item \cite{kato} If $r$ is a rational function with $r(\infty)=\infty$, and if $X = \{ z \in \mathbb C : |r(z)|\leq 1\}$
        is a convex set containing $W(A)$, then $w(r(A))\leq 1$.
      \end{enumerate}
 \end{enumerate}

\medskip

\noindent{\bf Open Problems:}
\begin{enumerate}
\item
  \cite{Cr07} Is $K_{Crouzeix}=2$? At least if we restrict ourselves to $3 \times 3$ matrices?
\end{enumerate}

\medskip

\ex
 \begin{enumerate}
 \item \cite{clark} The spectrum of the Toeplitz operator $T_\phi$ acting on the Hardy space $H^2$ with symbol $\phi(w)=aw + b w^{-1}$, $0<b<a$, is a convex set whose boundary is an ellipse with semi-axes $a\pm b$. The spectrum coincides with the closure of the numerical range $W(T_\phi)$. This set is not spectral but $K$-spectral for $T_\phi$ with  $K=\sqrt{1+|b/a|^2}$.

     Other classes of symbols with similar properties are given in \cite{clark}.
 \item  \label{c107s4ex0} \cite{PutinarSandberg}
     The configuration constant of Fact \ref{s5}.\ref{c17.s5.fa1b} can be estimated for
(smooth) $X$: $q(X)\leq 1-\frac{1+e}{2}\sqrt{1-e^2}$ for $\partial X$ an ellipse with eccentricity $e<1$ and
     $q(X) \leq 1 - \frac{\mbox{\small length}(\partial X)}{2\pi R}$, with $R$ the maximum radius of curvature.

   \item $K$-spectral sets containing the numerical range: let
$A\in \LH$, and let $X\subset \mathbb C$ be closed and convex containing $W(A)$.

 \begin{enumerate}
 \item \label{c107s5fa3}
    \cite{CD03, BCD06, BC07}
    If $X$ is a convex sector or a strip,
then $X$ is $(2+\frac{2}{\sqrt{3}})$-spectral for $A$.

 \item
    \cite{Cr04b, BC07}
    If the boundary of $X$ is a parabola or a hyperbola,
then $X$ is $(2+\frac{2}{\sqrt{3}})$-spectral for $A$.

 \item
    \cite[Theorem~1]{BC07}
    If the boundary of $X$ is an ellipse with eccentricity $e\leq 1$,
    then $X$ is $(2+\frac{2}{\sqrt{4-e^2}})$-spectral for $A$.

\item
   For an equilateral triangle, we have $K=2+\pi+6\log(2)$ according to Fact~\ref{s5}.\ref{c17.s5.fa1c},
   and for a square $K=2+\pi+4\log(2)$.

  \end{enumerate}

\item \cite{Cr13} For the $3\times 3$ matrix
$$A=\left[\begin{array}{c@{\quad}c@{\quad}c}
  0  &  2  &  0  \\
 \epsilon  &  0  &  (1-\epsilon^2)/\sqrt{2}  \\   0  &  0  &  (1-\epsilon^2)/\sqrt{2}  \end{array}\right]
$$
and sufficiently small $\epsilon > 0$, the numerical range $W(A)$ is $2$-spectral for $A$ but not completely $2$-spectral for $A$.
 \item
   \cite{greenbaumchoi} Crouzeix's conjecture is known to hold for generalized Jordan blocks where one
   replaces in a Jordan block of arbitrary size the lower left entry $0$ by an arbitrary scalar.
 \item
 \label{c107s4ex1}
   For $f(z)=z$ we recover from Fact \ref{s5}.\ref{c17.s5.fa2}
   the well-known link between numerical radius and euclidean norm, namely $[1,2] \ni \| A \|/w(A) = \| f(A) \|/\sup_{z\in W(A)} |f(z) |\leq K_{Crouzeix}$.
 \item
 \label{c107s4ex2} Consider  $A=\left[\begin{array}{c@{\quad}c}
0 & 2 \\0 & 0\end{array}\right]\in \mathbb C^{2\times 2}$,  $f(A)=\left[\begin{array}{c@{\quad}c}
f(0) & 2f'(0) \\0 & ~f(0)\end{array}\right]$. Here $W(A)$ is the closed unit disk and, for $f(z)=z$,
     $$
          \Bigl. \| f(A) \| \Bigr/ \sup_{z\in W(A)} |f(z) | = 2 \leq K_{Crouzeix}.
     $$
 \item
 \label{c107s4ex3}
 \cite[Section 2]{Cr04a} Consider $A=\left[\begin{array}{c@{\quad}c}1 & \rho-\frac{1}{\rho} \\0 & -1\end{array}\right]\in \mathbb C^{2\times 2}$ with $\rho>1$. 
 According to
 Fact 18.1.7 in Chapter 18,
 $W(A)$ is compact, with its boundary given by a ellipse with foci $\pm 1$
 and minor axis $\rho-1/\rho$.
The matrix can be diagonalized as
 $$
      A = L B L^{-1} , \quad B = \left[\begin{array}{c@{\quad}c}
1 & 0 \\0 & -1\end{array}\right] ,
      \quad L =\left[\begin{array}{c@{\quad}c}1 & -\frac{\rho-1/\rho}{\rho+1/\rho} \\0 & \frac{2}{\rho+1/\rho}\end{array}\right],
      \quad \| L \| \, \| L^{-1} \| = \| A \| = \rho .
 $$
 Thus
 $$
      \frac{\| f(A) \|}{\displaystyle \sup_{z\in W(A)} |f(z) |}
       \leq \rho \,
      \frac{\| f(B) \|}{\displaystyle \sup_{z\in W(A)} |f(z) |}
       = \rho \,
      \frac{\max\{ |f(-1)|, |f(1)|\}}{\displaystyle \sup_{z\in W(A)} |f(z) |} .
 $$
 By \cite[Theorem~2.1]{Cr04a} and some elementary calculus,
 the right-hand side of this expression is maximized for the function $f_0$ mapping conformally $W(A)$ on
the closed unit disk, with $f_0(0)=0$, $f_0'(0)>0$.
Thus $W(A)$ is a $K$-spectral set for $A$ with optimal constant $K=\rho \, f_0(1) = \| f_0(A) \|\in (1,2)$,
the last two relations following from the fact that an explicit formula for $f_0$ is known \cite[Eqn.\ (2.2)]{Cr04a}.
\end{enumerate}

\vspace*{1pc}
\section{Applications to the Approximate Computation\\ of Matrix Functions}\label{s6}

\smallskip

\fc

\smallskip

Notation: $A\in \mathbb C^{n \times n}$,
$\mathbb E\subset \mathbb C$ compact convex, being $K(\mathbb E)$--spectral for $A$, $f$ a function
being analytic on $\mathbb E$.

 \begin{enumerate}
 \item
 \label{c107s6fa1} Polynomial approximation through Taylor sums.

     A popular method \cite[Section~4.3]{Hi08} for approximately computing (entire) functions of
matrices is to approach the Taylor series $f(z)=\sum_{j=0}^\infty c_j z^j$ by its $m$th partial
sum $S_m(z)=\sum_{j=0}^m c_j z^j$, with the error estimate \cite[Cor.~2]{Ma93}
     $$
        \| (f-S_m)(A) \| \leq \frac{1}{(m+1)!} \max_{0\leq t \leq 1} \| A^{m+1} f^{(m+1)}(tA) \|.
     $$
     Here the right-hand side can by bounded in terms of $\mathbb E$ by the techniques of the
preceding subsections, see, e.g., Example~\ref{s1}.\ref{c107s5fa1} or Fact~\ref{s5}.\ref{c17.s5.fa2}.

 \smallskip

 \item
 \label{c107s6fa2} Polynomial approximation through Faber sums.

 \smallskip

     Instead of Taylor sums, following \cite[Section~4.4.1]{Hi08}
     and Section 11.7 one may also consider best
polynomial approximants $p_m$ of $f$ on
     $\mathbb E$, leading to the error estimate
     $$
        \| (f-p_m)(A) \| \leq K(\mathbb E) \, \rho_m(f,\mathbb E) , \quad
        \rho_m(f,\mathbb E) = \min_{\deg p \leq m} \max_{z\in \mathbb E} |f(z)-p(z)|.
     $$
     According to \cite[Theorem~4]{KP67}, this rate $\rho_m(f,\mathbb E)$ is achieved up to
some factor $\alpha\log(m)+\beta$ with explicit $\alpha,\beta>0$ not depending on $f$ nor on $\mathbb E$
     by taking as $p_m$ the $m$th partial Faber series, which is defined as follows: let $\phi=\psi^{-1}$
map conformally the exterior of $\mathbb E$ onto the exterior of the closed unit disk $\mathbb D$,
then the $m$th Faber polynomial $F_m$ is defined as the polynomial part of the Laurent
expansion of $\phi^m$ at $\infty$, and $f$ has the Faber series
     $$
          \sum_{j=0}^\infty f_j F_j(z) , \quad f_m=\frac{1}{2\pi i} \int_{|w|=1} \frac{f(\psi(w))}{w^{m+1}} dw,
     $$
     absolutely converging to $f$ uniformly in $\mathbb E$ \cite[Theorem~5]{KP67}.
     Notice that for $\mathbb E$ a disk centered at $0$ we recover Taylor series, and for $\mathbb E=[-1,1]$
Chebyshev orthogonal series \cite{Tr12}.
     By \cite[Section~3]{BR09} one has the a posteriori bound
     $|f_{m+1}|\leq \rho_m(f,\mathbb E) \leq 2(|f_{m+1}|+|f_{m+2}|+\ldots)$.
Estimates for the Faber coefficients of the exponential function can be found in \cite[Section~4]{BR09}.

     The Faber operator \cite{Ga87}
     $$
         \mathcal F (F)(z)=
         \sum_{j=0}^\infty f_j F_j(z), \quad F(w)=\frac{f_0}{2}+\sum_{j=1} f_j w^j
     $$
     for functions $F$ analytic on $\mathbb D$ maps polynomials of degree $m$ to
polynomials of degree $m$, in particular $\mathcal F(w^m)(z)=F_m(z)$ for $m \geq 1$. According to \cite{Be05},
     \cite[Theorem~2.1]{BR09}
     $$
          \| \mathcal F(F)(A) \| \leq 2 \, \max_{w\in \mathbb D} |F(w)| ,
     $$
     allowing to relate best polynomial approximation of $f$ on $\mathbb E$ to $\rho_m(\mathcal F^{-1}(f),\mathbb D)$,
and to derive error estimates for matrix functions which do not involve $K(\mathbb E)$, e.g., for the $m$th partial
Faber sum $p_m$ \cite[Theorem~3.2]{BR09}
     $$
          \| (f-p_m)(A) \| \leq 2 \, \sum_{j=m+1}^\infty |f_j| .
     $$

 \smallskip

 \item
 \label{c107s6fa3} Polynomial Arnoldi method.

 \smallskip

     The Arnoldi process is a popular method for approaching $f(A)b$
     for some fixed vector $b\neq 0$ and large sparse $A$ \cite[Section~13.2]{Hi08}. The $m$th approximant is given by
     the projection formula $V_m f(V_m^* A V_m) V_m^*b$, where $V_m\in \mathbb C^{n\times m}$ with columns
spanning an orthonormal basis of the Krylov subspace $\span\{ b,Ab,\ldots,A^{m-1}b\}$. Typically, $m \ll n$, and hence $f(V_m^* A V_m)$ can be computed by some direct method. According to \cite[Lemma~13.4]{Hi08}, this method is exact for $f$ a polynomial of degree $\leq m-1$, and hence
     $$
          \epsilon_m:= \frac{\| f(A)b - V_m f(V_m^* A V_m) V_m^*b \|}{\| b\|}
          \leq \| (f-p)(A) \| + \| (f-p)(V_m^* A V_m) \|
     $$
     for any polynomial $p$ of degree $\leq m-1$. Supposing that $W(A)\subset \mathbb E$
(and thus $W(V_m^* A V_m)\subset \mathbb E$), we obtain by means of the techniques of
     \ref{s6}.\ref{c107s6fa2} that $\epsilon_m \leq 2 \, K_{Crouzeix} \, \rho_{m-1}(f,\mathbb E)$  or
     $\epsilon_m \leq 4 \, \rho_{m-1}(\mathcal F^{-1}(f),\mathbb D) \leq 4 \sum_{j=m}^\infty |f_j|$,
see \cite[Proposition~3.1 and Theorem~3.2]{BR09}.

 \smallskip

 \item
 \label{c107s6fa4} Rational approximants with free poles via (Faber-)Pad\'e approximants.

 \smallskip

     Another popular approach (see \cite[Section~4.4.2]{Hi08} and Section 11.7) for approaching matrix functions for functions $f$
with singularities is to replace $f$ by the rational function $p/q$ with $p,q$ polynomials of degree at
most $k$, and $m$, respectively, $k \geq m-1$, such that the first $k+m+1$ terms in the
Taylor expansion of $f-p/q$ at zero vanish. Here one hopes that the
poles of this $[k|m]$ Pad\'e approximant $p/q$ do mimic the singularities of $f$.
Though this approach is also applied for entire functions like the exponential
function \cite[Section 10.7.4]{Hi08}, the error is best understood for Markov functions
     $$
         f(z) = c + \int_\alpha^\beta \frac{d\mu(x)}{z-x} , \quad
         c \in \mathbb R, \quad \mu \mbox{~some positive measure},
     $$
     where we suppose that $\beta<-w(A)$.
     This includes (up to some variable transformations) functions like $\log(z)$, $1/\sqrt{z}$
or more generally $p$th roots, $\mbox{sign}(z)$, $\tanh(z)$, and others \cite{BGM81,BR09}.
     By \cite[Lemma~6.2.1]{ST92}, the denominator $q$ has all its roots in $[\alpha,\beta]$,
and the error may be represented as
     $$
          f(z)-\frac{p}{q}(z) = \frac{z^{m+n+1}}{q(z)^2}
          \int_{\alpha}^\beta \frac{q(x)^2}{x^{m+n+1}} \frac{d\mu(x)}{z-x}.
     $$
     Then the error on the disk $\{ z \in \mathbb C : |z|\leq w(A)\}$ is minimal for $z=-w(A)$, and,
using Example \ref{s1}.\ref{c107s5fa1},
     $$
         \| (f-\frac{p}{q})(A) \| \leq 2 \, |(f-\frac{p}{q})(-w(A))|.
     $$

     In general, sharper error bounds are obtained for $k \geq m$ by combining the
above techniques with those of \ref{s6}.\ref{c107s6fa2}, where we suppose in addition
that $\mathbb E$ is symmetric with respect to the real axis. We first notice
that, with $f$, also $\mathcal F^{-1}(f)$ is a Markov function \cite[Theorem~6.1(a)]{BR09}.
Denoting by $P/Q$ the $[k|m]$ Pad\'e approximant of $\mathcal F^{-1}(f)$, the
function $\mathcal F(P/Q)$ is rational with numerator degree $\leq k$ and
denominator degree $\leq m$, called the $[k|m]$ Faber-Pad\'e approximant \cite{El83}.
For $\mathbb E=[-1,1]$ one recovers the so--called non-linear Chebyshev-Pad\'e approximant \cite{Su09}.
As above, we may bound for Markov functions the error through
     $$
         \| \Bigl(f-\mathcal F(\frac{P}{Q})\Bigr)(A) \| \leq 2 \,
         | \Bigl(\mathcal F^{-1}(f)-\frac{P}{Q}\Bigr)(-1) |.
     $$

 \smallskip

 \item
 \label{c107s6fa5} Rational approximation with prescribed poles and rational Arnoldi.

 \smallskip

     There exists a variant of the polynomial Arnoldi method \ref{s6}.\ref{c107s6fa3} where
the columns of $Q_m\in \mathbb C^{n \times m}$ give an orthonormal basis of the rational
Krylov subspace $$q(A)^{-1}\span\{ b,Ab,\ldots,A^{m-1}b\}$$ for some fixed polynomial $q(z)=\prod_j (z-z_j)$
of degree $\leq m-1$ (and hence for $q=1$ we recover the polynomial Arnoldi method).
The computation of $V_m$ and $V_m^* A V_m$ by some rational variant of the Arnoldi process \cite{BR09}
requires $(A-z_j I)^{-1}a$ for some vectors $a$, and this task of solving shifted linear systems is
particularly trackable if $q$ has a small number of multiple poles.
     As before we obtain
     $$
          \widetilde \epsilon_m :=\frac{\| f(A)b - V_m f(V_m^* A V_m) V_m^*b \|}{\| b \|}
          \leq \| (f-\frac{p}{q})(A) \| + \| (f-\frac{p}{q})(V_m^* A V_m) \|
     $$
     for any polynomial $p$ of degree $\leq m-1$ \cite[Theorem~5.2]{BR09}. Let $W(A)\subset \mathbb E$,
and suppose that $\mathbb E$ is symmetric with respect to the real axis, and $z_j \not\in \mathbb E$.
For estimating the error we are left with the task of approaching $f$ on $\mathbb E$ by a rational
function with fixed denominator $p/q$.
     Notice that $p/q = \mathcal F(P/Q)$ with $P$ a polynomial of degree at most $\leq m-1$,
and $Q(w)=\prod_j (w-\phi(z_j))$ \cite{El83}.
     Hence $\widetilde \epsilon_m\leq 4 \, \rho_{m-1}^Q(\mathcal F^{-1}(f),\mathbb D)$, where
     $$
          \rho_{m-1}^q(f,\mathbb E)= \min_{\deg p \leq m-1} \max_{z\in \mathbb E} |f(z)-\frac{p}{q}(z)|,
     $$
     see \cite[Theorem~5.2]{BR09}. For Markov functions $f$ as in
     \ref{s6}.\ref{c107s6fa4}, lower and upper bounds for $\rho_{m-1}^Q(\mathcal F^{-1}(f),\mathbb D)$
are given in \cite[Theorem~6.2]{BR09}, in particular the simple explicit bound
     $$
        \rho_{m-1}^Q(\mathcal F^{-1}(f),\mathbb D) \leq
        \frac{1}{|\phi(\beta)|} \, \, \max_{z\in \mathbb E} | f(z) - f(\infty) |    \, \,
\max_{w\in [\phi(\alpha),\phi(\beta)]}
        \Bigl| \prod_j \frac{w-\phi(z_j)}{1-w\overline{\phi(z_j)}} \Bigr|.
     $$
     Poles $z_j$ minimizing the right-hand should therefore be in $[\alpha,\beta]$,
the set of singularities of $f$, and various configurations of poles minimizing the
right-hand side have been considered in \cite[Section~6]{BR09}.

 \smallskip

 \item
 \label{c107s6fa6} Error bounds for GMRES and FOM for solving systems of linear equations.

 \smallskip

    Both GMRES and FOM are iterative Krylov subspace methods for solving systems $A x = b$ with $A$
large and sparse, see Chapter 41. Here we may apply the techniques of the preceding sections for the Markov function $f(z)=1/z$ provided that $0\not\in \mathbb E$, $\mathbb E$ containing $W(A)$.
    The residual of the $m$th iterate $x_m^{GMRES}$ of GMRES with starting residual $r=b-Ax_0^{GMRES}$ satisfy \cite{Be05}
    $$
         \frac{\| b - Ax_m^{GMRES}\|}{\|r\|} = \min_{\deg p \leq m} \frac{\| p(A)r\|}{|p(0)|\, \|r\|}
         \leq \min\{ 1 , \frac{2}{|F_m(0)|} \} \leq \frac{2+1/|\phi(0)|}{|\phi(0)|^m}
    $$
    with $F_m$ the $m$th Faber polynomial of $\mathbb E$, and $\phi$ mapping conformally the exterior
of $\mathbb E$ onto the exterior of the unit disk. The asymptotic convergence factor $1/|\phi(0)|<1$
can be computed for various shapes of $\mathbb E$. For instance \cite{Be05,BGT06}, for positive
definite $A+A^*$ considering the lens $\mathbb E=\{ z\in \mathbb C : \mathrm{re}(z)\geq \mbox{dist}(0,W(A)), |z|\leq w(A) \}$
we get $1/|\phi(0)|=2\sin(\beta/(4-2\beta/\pi))<\sin(\beta)$,
with the angle
    $\beta\in (0,\pi/2)$ being defined by $\cos(\beta) = \frac{\mbox{dist}(0,W(A))}{w(A)}$.

    The $m$th iterate of FOM (with starting vector $x_0^{FOM}$ is a special case of the polynomial Arnoldi method
    \ref{s6}.\ref{c107s6fa3}, namely
    $x_m^{FOM}=V_m (V_m^* A V_m)^{-1} V_m^* b$ or $f(z)=1/z$, and thus
    $$
        \frac{\| x_m^{FOM} - A^{-1} b\|}{\| b \|} \leq 4 \eta_{m-1}(\mathcal F^{-1} (f),\mathbb D) \leq
       \frac{4 |\phi(0)|^{-m}}{\mbox{dist}(0,\mathbb E)}.
    $$

\end{enumerate}

\smallskip



\end{document}